\documentclass[12pt]{article}
\usepackage{a4,amsmath,amssymb,cite,ytableau}
\usepackage{rotating}
\usepackage{tikz}
\newtheorem{theorem}{Theorem}[section]
\newtheorem{mydef}[theorem]{Definition}
\newtheorem{prop}[theorem]{Proposition}
\newtheorem{conj}[theorem]{Conjecture}
\def\Bbb{\mathbb}
\def\BZ{\Bbb Z}

\catcode`@=11 \@addtoreset{equation}{section} \catcode`@=12

\begin{document}
\bibliographystyle{utphys}
\begin{titlepage}
\renewcommand{\thefootnote}{\fnsymbol{footnote}}
\noindent
{\tt IITM/PH/TH/2011/7}\hfill
{\tt arXiv:1203.nnnn} \\[4pt]
\mbox{}\hfill 
\hfill{\fbox{\textbf{v1.0; March 2012 }}}

%\vspace{1.0cm}
\begin{center}
\large{\sf  Notes on  higher-dimensional partitions}
\end{center} 
\bigskip 
\begin{center}
{\sf Suresh Govindarajan\footnote{suresh@physics.iitm.ac.in} } \\[3pt]
\textit{Department of Physics, Indian Institute of Technology Madras,\\ Chennai 600036, India \\[4pt]
}
\end{center}
\bigskip
\bigskip
\begin{abstract}
We show the existence of a series of transforms that capture several 
structures that underlie  higher-dimensional partitions. These transforms 
lead to a sequence of triangles whose entries are given combinatorial 
interpretations as the number of particular types of skew Ferrers 
diagrams. The end result of our analysis is the existence of a triangle, that
we denote by $F$, which implies that the data needed to compute
the number of partitions of a given positive integer is reduced by a factor
of half.  The number of spanning rooted forests appears intriguingly in 
a family of entries in the $F$.
Using modifications of an algorithm due to Bratley-McKay, we are
able to directly enumerate entries in some of the triangles. As a result, we
have been able to compute numbers of partitions of positive integers $\leq 25$ in any dimension. 
\end{abstract}
\end{titlepage}
\setcounter{footnote}{0}
 \newcommand{\nodee}[4]{\begin{smallmatrix} #1\\ #2\\  #3 \\ #4 \end{smallmatrix}}
\section{Introduction} 

An unrestricted $d$-dimensional partition of $n$ is a collection of $n$ 
points (nodes) in $\mathbb{Z}_+^{d+1}$ satisfying the following 
property: if the collection contains a node $\mathbf{a}=(a_1,a_2,\ldots, 
a_{d+1})$, then all nodes $\mathbf{x}=(x_1,x_2,\ldots,x_{d+1})$ with 
$0\leq x_i \leq a_i$ $\forall$ $i=1,\ldots,d+1$ also belong to 
it\cite{Atkin:1967,AndrewsPartitions}. Let $p_d(n)$ denote the number of 
distinct such partitions. Denote by $P_d(q)$, the generating function of 
unrestricted $d$-dimensional partitions: ($p_d(0)\equiv 1$)
\begin{equation}
P_d(q)=\sum_{n=0}^\infty \ p_d(n)\ q^n \ .
\end{equation}
There exist explicit formulae for the generating functions for $d=1$ and 
$d=2$ due to Euler and MacMahon respectively\cite{MacMahon}.  However, 
no such formulae exist for $d>2$ as an inspired guess of MacMahon was 
subsequently proven to be false\cite{Atkin:1967}. It appears that there 
is no simple formula and one has to take recourse to brute force 
enumeration. Given that asymptotically one 
has\cite{Wu:1996,Bhatia:1997,Balakrishnan:2011bm}
\begin{equation}
\log p_d(n) \sim n^{d/d+1}\ ,
\end{equation}
it is easy to see that the numbers grow exponentially fast and naive 
enumeration is not the way to go.

The first serious attempt at direct enumeration of partitions in any 
dimension is due to Atkin et. al. \cite{Atkin:1967} based on an 
algorithm due to Bratley and McKay\cite{Bratley:1967a}. Knuth provided 
another algorithm that enumerates numbers of topological sequences which 
can be used, in principle, to generate numbers of partitions in any 
dimension\cite{Knuth:1970}. Both algorithms are highly recursive and 
easily implemented on a computer.

This paper attempts to find structures in the enumeration of partitions 
and come up with refinements in their enumeration. Such refinements when 
cleverly combined with computer-based enumeration should in principle 
enable one to enumerate partitions of integers below some maximum value 
in any dimension. The maximum value turns out be $25$ in our case though 
we believe that, with some effort, this number can be pushed to around 
$30$.

Our refinements begin with the result of of Atkin et. al. who showed 
that the binomial transform of $p_d(n)$ leads to a lower-triangular 
matrix that we denote by $A=(a_{n,r})$.
\begin{equation*}
p_d(n) = \sum_{r=0}^{d+1} \binom{d+1}{r} \ a_{n, r} \ .\eqno{\eqref{Amatrixdef}}
\end{equation*}
This transform implies that in order to compute partitions of a positive 
integer $n$ in any dimension, we need to only compute $(n-1)$ numbers 
that make up a particular row of the triangle $A$. 
We show the existence another triangular matrix, that we denote by 
$F=(f_{n,x})$, as a transform of the matrix $A$ with fewer entries.
\begin{equation*}
a_{m+r+1,r} = \sum_{x=0}^r \sum_{p=x}^{m} \binom{r}{x} \binom{\tbinom{r-x}2}{m-p} \ f_{p+x+1,\ x}\ .
\eqno{\eqref{Fdef}}
\end{equation*}
Our result is that we need only $[(n-1)/2]$ independent numbers i.e., 
roughly half of the initial estimate to determine partitions of $n$ in 
any dimension. We illustrate the gain by explicitly displaying the first eleven rows of the $A$ and $F$-matrices.
$$
A=\left(
\begin{smallmatrix}
 1  \\
 0 & 1    \\
 0 & 1 & 1   \\
 0 & 1 & 3 & 1   \\
 0 & 1 & 5 & 6 & 1   \\
 0 & 1 & 9 & 18 & 10 & 1 \\
 0 & 1 & 13 & 44 & 49 & 15 & 1 \\
 0 & 1 & 20 & 97 & 172 & 110 & 21 & 1  \\
 0 & 1 & 28 & 195 & 512 & 550 & 216 & 28 & 1  \\
 0 & 1 & 40 & 377 & 1370 & 2195 & 1486 & 385 & 36 & 1 \\
 0 & 1 & 54 & 694 & 3396 & 7603 & 7886 & 3514 & 638 & 45 & 1  \\
\end{smallmatrix}
\right)\quad,\quad
F=\left(\begin{smallmatrix}
1\\
0  \\
 0&1  \\
0& 1  \\
 0&1 & 3  \\
 0&1 & 7  \\
 0&1 & 11 & 16  \\
 0&1 & 18 & 58  \\
 0& 1 & 26 & 135 & 125  \\
 0& 1 & 38 & 293 & 618 \\
 0& 1 & 52 & 574 & 1927 & 1296 
   \end{smallmatrix}\right)
\ .
$$

The $F$-triangle is, in a sense, the end-point of a sequence of transforms 
and triangles that we introduce. We also provide combinatorial 
interpretations for the various triangles that appear as a result of 
these transforms. This enables use to modify the Bratley-McKay(BM) 
algorithm to directly enumerate the matrix $A$ that we mentioned earlier 
and a second triangle, $C$ that we define in the sequel. As we discuss in the appendix, similar refinements can be carried out for partitions restricted in a box.

\subsection{Summary of results}
\begin{enumerate}
\item Given a partition in any dimension, we have introduced two new 
attributes: its intrinsic dimension (i.d.) - see definition \ref{iddef} 
and its reduced dimension (r.d.) - see definition \ref{rddef}.
\item These two attributes lead to two new triangles, the $A$ and 
$C$-matrices (see Eq. \eqref{Amatrixdef} and \eqref{Cmatrixdef}) whose 
entries admit combinatorial interpretations. We propose a further 
refinement in the form of two other triangles, the $D$ matrix(see Eq. 
\eqref{Dmatrixdef}).
\item We show that the $C/D$ triangles are the first in a series of transforms, the end-point of which 
leads to a triangle $F$ (see Eq. \eqref{Fdef}). The $n$-th row of this matrix has only $[(n-1)/2]$  entries (where $[x]$ is the integral part of $x$) and these entries determine the partitions of $n$ in any dimension. This constitutes the main result of this paper.
\item We see an intriguing  relationship between the numbers of spanning rooted forests on $m$ vertices and $\alpha$ components and a family of entries in the $F$-matrix. This is Proposition \ref{Cayley}.
\item We conjecture the existence of two other triangles, the $\alpha$- 
and the $\beta$-matrices with integer entries.
\item We prove a conjecture of Hanna on the existence of a triangle that 
determines all higher-dimensional partitions.
\item We propose a modification to an algorithm of Bratley and McKay 
that enables us to directly compute the $A$ and $C$ matrices. We compute 
the first 25 rows of the $F$-matrix 
thereby obtaining partitions in all dimensions for integers $\leq 25$. 
\item Tables \ref{Knuthresults}-\ref{Fboxtable} provide the numerical 
results that we have obtained.
\end{enumerate}

\section{Structures in higher-dimensional partitions}

\subsection{Ferrers diagrams and permutation symmetry}

A Ferrers diagram represents the partition as a $(d+1)$-dimensional 
arrangement of nodes. For instance, the following one-dimensional 
partition of $4$
\begin{displaymath}
\left\{
\begin{pmatrix}
 0 \\ 0 
\end{pmatrix}, \ 
 \begin{pmatrix}
 1 \\ 0 
\end{pmatrix}, \ 
\begin{pmatrix}
 0 \\ 1 
\end{pmatrix}, \ 
\begin{pmatrix}
 0 \\ 2
\end{pmatrix}
\right\}  \textrm{ or }
\left(\begin{smallmatrix}0&1&0&0 \\ 0 & 0 & 1 & 2\end{smallmatrix}\right) \textrm{ in compressed form}\ ,
\end{displaymath}
is represented by the following two-dimensional Ferrers diagram or as a 
Young diagram where we replace the nodes by squares(more generally, 
hypercubes).

\begin{center}
\begin{tikzpicture}[scale=0.4]
% Draw x and y axis lines
\draw [->] (-0.3,0.3) -- (2.8,0.3) node [right] {$x_2$};
\draw [->] (-0.3,0.3) -- (-0.3,-1.8) node [below] {$x_1$};
\node  at (0,0) [fill,circle,inner sep=1.3pt] {};
\node at (0,-1) [fill,circle,inner sep=1.3pt] {};
\node at (1,0) [fill,circle,inner sep=1.3pt] {};
\node at (2,0) [fill,circle,inner sep=1.3pt] {};
\node at (8,-1) {or};
\end{tikzpicture}
%\vspace{1cm}
\ytableausetup{aligntableaux=bottom} \qquad \ydiagram{3,1}
\end{center}

There is a natural action of $S_{d+1}$ on the $(d+1)$-dimensional 
Ferrers diagram -- this corresponds to permuting the $(d+1)$ 
coordinates. For one-dimensional partitions, this is referred to as 
conjugation. The symmetry group of a $d$-dimensional partition is the 
\textit{largest} sub-group of $S_{d+1}$ that acts trivially on the 
corresponding Ferrers diagram.

\subsection{The intrinsic dimension}

Typically, one is interested in the asymptotic behavior of $p_d(n)$ for 
large number of nodes $n$ while keeping the dimension $d$ fixed. 
However, one may ask about what happens to $p_d(n)$ if we keep the 
number of nodes. i.e., $n$, fixed and keep increasing $d$. It is easy to 
see that when $d>n+1$, all the nodes of the Ferrers diagram necessarily 
lie in some $r$-dimensional hyperplane with $r<d$. This motivates the 
following definition (implicitly present in Atkin et. 
al.\cite{Atkin:1967}).

\begin{mydef}\label{iddef}
Given a Ferrers diagram, let it be contained in a $r$-dimensional 
hyperplane but not in any $(r-1)$-dimensional hyperplane. The intrinsic 
dimension(i.d.) of the Ferrers diagram is defined to be $r$.
\end{mydef}

Note that such a $r$-dimensional hyperplane is given by setting 
$(d+1-r)$ coordinates to zero. Any permutation of the $(d+1-r)$ 
coordinates (that are set to zero to obtain the hyperplane containing 
the nodes) does not change the Ferrers diagram. It is thus easy to see 
that the symmetry of a Ferrers diagram in $(d+1)$-dimensions of i.d. $r$ 
is necessarily of the form $H\times S_{d+1-r}$ where $H\subseteq S_{r}$. 
We shall (somewhat loosely) call $H$, the symmetry of the Ferrers 
diagram.

Let two $d$-dimensional partitions be equivalent if their Ferrers 
diagram are related by an $S_{d+1}$ action. It is easy to see that all 
$d$-dimensional partitions belonging to such an equivalence class have 
the same intrinsic dimension. Further, given a $d+1$-dimensional Ferrers 
diagram with symmetry $H$ and i.d. $r$, the number of Ferrers diagrams 
in its equivalence class is given by the order of the coset 
$S_{d+1}/(H\times S_{d+1-r})$ i.e.,
$$
 \frac{(d+1)!}{(d+1-r)!\times \textrm{ord}(H)} = \binom{d+1}{r} \times \frac{r!}{\textrm{ord}(H)}\ .
$$
 \begin{mydef}
A Ferrers diagram is said to be strict when its intrinsic dimension 
equals its dimension.
 \end{mydef}
Given a $d+1$-dimensional Ferrers diagram of i.d. $r$, it is useful to 
drop the $(d+1-r)$ dimensions that are orthogonal to the hyperplane 
containing the nodes thus obtaining a strict FD.  The symmetry of the 
strict Ferrers diagram is now $H\subseteq S_r$.
\begin{mydef} 
A generalized Ferrers diagram (gFD) refers to the equivalence class of 
strict Ferrers diagrams obtained by the action of $S_r$ on a given 
strict Ferrers diagram of i.d. $r$.
\end{mydef}
The number of strict FD's in a gFD of i.d. $r$ and symmetry group $H$ is 
$\tfrac{r!}{\textrm{ord}(H)}$.
\begin{mydef}
The weight of a gFD of i.d. $r$ and symmetry $H\subseteq S_r$ is defined 
to be $\tfrac{r!}{\textrm{ord}(H)}$.
\end{mydef}
Since $H\subseteq S_r$, Lagrange's theorem implies that the weight, 
$\frac{r!}{\textrm{ord}(H)}$, is a positive non-zero integer.  Note that 
the weight is independent of the dimension of the Ferrers diagram and is 
the same for \text{all} elements in an equivalence class. Thus, to an 
equivalence class of a given Ferrers diagram, we associate three 
numbers: the number of nodes $n$, the i.d. $r$, and the weight, $w$.  An 
important observation is that there exist \textit{no} Ferrers diagram 
with $n$ nodes and i.d. $r\geq n$ -- this follows from noting that one 
needs at least $r+1$ nodes to create a Ferrers diagram of i.d. $r$. We 
see that the number of $d$-dimensional partitions is thus given by
\begin{align}
p_d(n) &=  \sum_{r=0}^{n-1}  \binom{d+1}{r} \sum_{\lambda\vdash (n,r)} 1 \\ 
&=\sum_{r=0}^{n-1}  \binom{d+1}{r} \sum_{[\lambda]\vdash (n,r)}  w(\lambda) \\ &:= \sum_{r=0}^{n-1} \binom{d+1}{r} a_{n,r}\ ,
\end{align}
where the second line defines $a_{n,r}$ as the sum over all strict FD's 
with $n$ nodes and i.d. $r$. In the second line, the sum over 
$[\lambda]$ indicates that we sum over equivalence classes of strict 
Ferrers diagrams(gFD). Note that $a_{n,r}$ has no dependence on $d$ and 
counts the numbers of strict Ferrers diagrams with $n$ nodes and i.d. 
$r$. We shall provide a second, and more useful, combinatorial 
description of $a_{n,r}$ later.

\subsection{The first transform}

We extend $a_{nr}$ into a lower-triangular matrix, that we denote by 
$A$, by setting $a_{nr}=0$ when $r\geq n$. Thus, we obtain the matrix 
$A=\big(a_{nr}\big)$ for $n=1,2,\ldots$ and $r=0,1,2,\ldots$. With this 
definition, we an rewrite the above equation as
\begin{equation}
\boxed{
p_d(n) = \sum_{r=0}^{d+1} \binom{d+1}{r} a_{nr}
}\ . \label{Amatrixdef}
\end{equation}
To our knowledge, the above observation first appeared in a paper by 
Atkin et. al.\cite{Atkin:1967}. Thus the $p_d(n)$, for a fixed value of 
$n$, corresponds to the Binomial Transform of the $n$-th row of the 
matrix $A$. It is easy to see that $a_{n,0}=\delta_{n,1}$. The lower 
triangular nature of $A$ implies that only $(n-1)$ numbers, 
$(a_{n,1},a_{n,2},\ldots, a_{n,(n-1)})$ determine $p_d(n)$ for any $d$. 
The matrix A appears in the OEIS as sequence number A119271\cite{oeis}. 
The inverse Binomial transform is given by
\begin{equation}
\boxed{
a_{nr} =\sum_{d=0}^{r-1} (-1)^{d+r+1} \binom{r}{d+1}\ p_d(n)\quad \textrm{for } n\geq r+1} \ ,
\end{equation}
with $ p_0(n)\equiv1 $. Of course, $ a_{nr}=0 $ when $ n<r+1 $ 
reflecting the lower-triangular nature of the matrix. Suppose we know 
all partitions of $ n_{max} $ up to $ d_{max} $. This determines the 
first $ n_{max}$ rows and $(d_{max}+1)$ columns of the matrix $A$.

For low values of $n$, we can explicitly compute the entries in the 
$A$-matrix by listing the gFD's and working out their weights as we do 
below.
\begin{align}
\ytableausetup{smalltableaux}
p_d(2)&=\tbinom{d+1}1 w\Big(\ydiagram{2}\Big) %\nonumber \\
=\tbinom{d+1}1  \\
p_d(3) &= \tbinom{d+1}1 w\Big(\ydiagram{3}\Big)+\tbinom{d+1}2 w\bigg(\ydiagram{2,1}\bigg) 
%\nonumber \\
=\tbinom{d+1}1 + \tbinom{d+1}2\ . \\
p_d(4) &= \tbinom{d+1}1 w\Big(\ydiagram{4}\Big)+
\ytableausetup{aligntableaux=center}\tbinom{d+1}2 w\bigg(\ydiagram{3,1}\bigg) +\tbinom{d+1}2 w\bigg(\ydiagram{2,2}\bigg) + \tbinom{d+1}3 
w\bigg(\!\! \includegraphics[height=7mm]{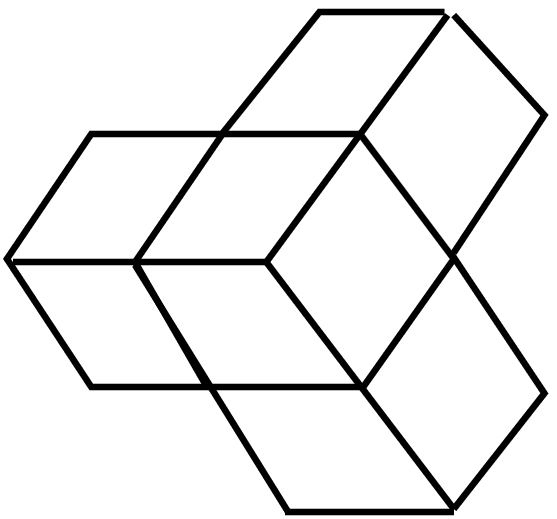}\bigg) 
\nonumber \\
&=\tbinom{d+1}1 + 3\tbinom{d+1}2+\tbinom{d+1}3\ .
\end{align}
The first few rows of the $A$-matrix are as follows (see also 
\cite{Ekhad:2012})
$$
A=\left(
\begin{array}{ccccccccccccc}
 1  \\
 0 & 1    \\
 0 & 1 & 1   \\
 0 & 1 & 3 & 1   \\
 0 & 1 & 5 & 6 & 1   \\
 0 & 1 & 9 & 18 & 10 & 1 \\
 0 & 1 & 13 & 44 & 49 & 15 & 1 \\
 0 & 1 & 20 & 97 & 172 & 110 & 21 & 1  \\
 0 & 1 & 28 & 195 & 512 & 550 & 216 & 28 & 1  \\
 0 & 1 & 40 & 377 & 1370 & 2195 & 1486 & 385 & 36 & 1 \\
 0 & 1 & 54 & 694 & 3396 & 7603 & 7886 & 3514 & 638 & 45 & 1  \\
 0 & 1 & 75 & 1251 & 7968 & 23860 & 35115 & 24318 & 7484 & 999 & 55 & 1   
\end{array}
\right)\ .
$$
 \begin{mydef}
Consider a pair of FD's $(\lambda,\mu)$ such that $\mu\subseteq 
\lambda$. Then, a skew Ferrers diagram is the set of nodes $\lambda 
\setminus \mu$.
\end{mydef}
One can think of the entries in the $A$-matrix as counting skew Ferrers 
diagrams obtained by deleting the node at the origin $(0,0,\ldots,0)^T$ 
that is contained in any Ferrers diagram. Then, $a_{n,r}$ is the number 
of strict FD's of dimension $r$ obtained by adding $(n-1)$ nodes to the 
node at the origin.  One sets $a_{1,0}\equiv1$.

 \subsection{A combinatorial interpretation}
 
We will now provide another combinatorial interpretation for the numbers 
$a_{nr}$ that make up the lower-triangular matrix $A$.  We begin with 
the observation that $a_{r+1,r}=1$ -- this follows because there is a 
unique FD of i.d. $r$ containing $r+1$ nodes. The coordinates are given 
in the following $r\times (r+1)$ matrix\footnote{Recall that each column 
is the coordinate of a node and thus there are $(r+1)$ columns and $r$ 
rows.}
\begin{equation}\label{specialFD}
 \mu_r := \left(\nodee00\vdots0 \nodee 10\vdots0 \nodee 01\vdots0 \ldots \nodee 00\vdots 1\right)\ .
 \end{equation}
 This FD has maximal symmetry $S_r$ and weight $1$. \\
 
\noindent\textbf{Remark:} Every FD with intrinsic dimension $r$ 
necessarily contains $\mu_r$. This implies that an FD with $n$ nodes and 
i.d. $r$ can be obtained by adding $m=n-r-1$ additional nodes to 
$\mu_r$. This leads to the following combinatorial interpretation for 
$a_{m+r+1,r}$.
 \begin{prop}
$a_{m+r+1,r}$ is the number of strict Ferrers diagrams with i.d. $r$ 
obtained by adding $m$ nodes to the standard Ferrers diagram, $\mu_r$.
 \end{prop}
Let $\lambda$ be an FD that contributes to $a_{n,r}$. Its symmetry group 
$H\subseteq S_r$ -- this implies that there will be 
$r!/\textrm{odd}(H)=wt(\lambda)$ distinct FD's obtained from it by the 
action of $S_r$. It is easy to see that the process of adding $m$ nodes 
to $\mu_r$ will generate the same number of FD's that belong to the 
equivalence class (gFD) $[\lambda]$.
 
So far we have completely determined the first 25 rows of the $A$-matrix 
(see Table \ref{Amatrix}). The entries have been determined by combining 
several methods: (i) taking the inverse Binomial transform of known 
numbers for higher-dimensional partitions, (ii) by direct enumeration 
using the combinatorial interpretation and (iii) by determining another 
triangle, the $C$-matrix, that we introduce later.  It is important to 
note that the numbers, when available, from the different methods agree. 
Further, none of the conjectural formulae are used in determining the 
entries.

 \subsection{The second transform}

\begin{mydef}\label{rddef}
Let $\lambda$ be an FD of i.d. $r$ and consider the skew FD 
$\lambda\setminus \mu_r$. Let the nodes of the skew FD be contained in a 
$x$-dimensional hyperplane (obtained by setting $r-x$ coordinates to 
zero) but not in any $(x-1)$-dimensional hyperplane. The reduced 
dimension (r.d.) of the FD $\lambda$ is said to be $x$.
 \end{mydef}
Clearly the reduced dimension of an FD is always less than or equal to 
its intrinsic dimension.  The symmetry of a FD with i.d. $r$ and r.d. 
$x$ is necessarily of the form $H\times S_{r-x}\subset S_r$. Then, one 
has
\begin{equation}
 a_{m+r+1,r}= \sum_{x=0}^{r} \binom{r}{x} c_{m,x}\ ,
 \end{equation}
where the binomial term $\binom{r}{x}$ takes into account the situation 
with maximal symmetry and $c_{0,0}\equiv 1$ and $c_{m,0}=c_{0,m}\equiv0$ 
for $m>0$.
\begin{enumerate}
\item The coefficients $c_{m,x}$ are clearly independent of the the i.d. 
($r$) as they are related to the skew FD's with $m$ nodes and r.d. $x$.
\item We say that a skew FD is \textit{strict} if its dimension and r.d. 
are the same.
\item Let us denote the equivalence class of strict skew Ferrers diagrams, $\lambda \setminus \mu_x$, under the $S_x$ action as an sFD.  All skew FD's  in an sFD  will have identical reduced and intrinsic dimensions.
Thus, given such a skew Ferrers diagram with 
symmetry $H\subseteq S_x$, its equivalence class will contain 
$\tfrac{x!}{\textrm{ord}(H)}$ distinct skew Ferrrers diagrams.
\item The $c_{m,x}$ are non-negative integers since they count the 
number of strict skew FD's with $m$ nodes and r.d. $x$.
\item For fixed $m$, one can see that the maximum value of r.d. with $m$ 
nodes is $2m$.  This enables us to convert the above equation into a 
second binomial transform
 \begin{equation}
 \boxed{
 a_{m+r+1,r}= \sum_{x=0}^{2m} \binom{r}{x} c_{m,x}}\ , \label{Cmatrixdef}
 \end{equation}
where we extend $c_{m,x}$ into a triangle, $C=(c_{m,x})$, by setting 
$c_{m,x}=0$ for $x>2m$.  We usually do not write out the zeroth row and 
column of the $C$-matrix.
\item For fixed $m$, we can consider $a_{m+r+1,r}$ as a function of 
$r$. The function $g_m(r):= 2m!!\ a_{m+r+1,r}$ is a polynomial of degree 
$2m$, conjecturally with integer coefficients, in the variable $r$ and 
$g_m(0)=0$ for $m>0$.
 \item We have directly determined eleven rows ($m\in[0,10]$) of the $C$-matrix (see Table 
\ref{Cmatrix}). The first few rows of the $C$-matrix are:

 $$
C= \begin{pmatrix}
1 \\
0&  1 & 1  \\
0& 1 & 3 & 6 & 3  \\
 0& 1 & 7 & 20 & 46 & 45 & 15  \\
0&  1 & 11 & 61 & 198 & 480 & 645 & 420 & 105  \\
 0& 1 & 18 & 138 & 706 & 2508 & 6441 & 10395 & 9660 & 4725 & 945  
\end{pmatrix}\ .
 $$
\end{enumerate}
It is easy to see that there is only one sFD with $m$ nodes and r.d. 
$2m$. In the picture below, the $m$ nodes of the sFD are indicated by 
open circles. The filled circles indicate the nodes of $\mu_{2m}$ that 
must be added to the sFD to obtain an FD.
\begin{equation}
\begin{gathered}
\begin{tikzpicture}[scale=0.3,rotate=270]
\draw [->] (0,1) -- (-2,3) node [right] {$x_2$};
\draw [->] (0,1) -- (2,3) node [below] {$x_1$};
\node  at (0,1) [fill,circle,inner sep=1.3pt] {};
\node  at (0,3)  {$\circ$};
\node at (1,2) [fill,circle,inner sep=1.3pt] {};
\node at (-1,2) [fill,circle,inner sep=1.3pt] {};
\node at (0,7) {$\times\quad$};
\end{tikzpicture}
\begin{tikzpicture}[scale=0.3,rotate=270]
\draw [->] (0,1) -- (-2,3) node [right] {$x_4$};
\draw [->] (0,1) -- (2,3) node [below] {$x_3$};
\node  at (0,1) [fill,circle,inner sep=1.3pt] {};
\node  at (0,3)  {$\circ$};
\node at (1,2) [fill,circle,inner sep=1.3pt] {};
\node at (-1,2) [fill,circle,inner sep=1.3pt] {};
\node at (0,9) {$\quad\times\quad \cdots\quad \times $};
%\node at (0,8) {};
\end{tikzpicture}
\begin{tikzpicture}[scale=0.3,rotate=270]
\draw [->] (0,1) -- (-2,3) node [right] {$x_{2m}$};
\draw [->] (0,1) -- (2,3) node [below] {$x_{2m-1}$};
\node  at (0,1) [fill,circle,inner sep=1.3pt] {};
\node  at (0,3)  {$\circ$};
\node at (1,2) [fill,circle,inner sep=1.3pt] {};
\node at (-1,2) [fill,circle,inner sep=1.3pt] {};
\end{tikzpicture}
\end{gathered}\label{reduciblesFD}
\end{equation}
The symmetry of the skew FD is $(S_m \ltimes \BZ_2^m)$ and thus 
$c_{m,2m}$ is the dimension of the coset i.e.,
$$
c_{m,2m}=\frac{\textrm{dim}(S_{2m})}{\textrm{ord}(S_m \ltimes \BZ_2^m))}=\frac{2m!}{2m!!}=(2m-1)!!\ .
$$

\begin{mydef} A skew FD of i.d. $r$ is said to be reducible if a proper 
subset of its nodes are contained in a $d$-dimensional hyperplane 
(obtained by setting $r-d$ coordinates to zero) with $d<r$ and the nodes 
not in the proper subset lie in the orthogonal $(r-d)$-dimensional 
hyperplane (obtained by setting the other $d$ coordinates to zero). 
\end{mydef} 
\begin{mydef} We say that an FD, $\lambda$, of 
i.d. $r$ is reducible if the skew FD, $\lambda\setminus\mu_r$ is 
reducible.
\end{mydef}
Thus a reducible sFD has multiple \textit{components} consisting of 
non-intersecting proper subsets of its nodes lying in mutually 
orthogonal hyperplanes. Thus the sFD given in Eq. \eqref{reduciblesFD} 
is reducible with $m$ components each of which is isomorphic to the 
irreducible sFD $\sigma_2$ defined as follows:
\begin{equation}
\begin{gathered}
\begin{tikzpicture}[scale=0.2]
\node at (-3,0.1) {$\sigma_{2}\equiv$};
\draw [->] (1,0) -- (3,-2) node [right] {};
\draw [->] (1,0) -- (3,2) node [below] {};
\node  at (1,0) [fill,circle,inner sep=1.3pt] {};
\node  at (3,0)  {$\circ$};
\node at (8,0){$=\left(\begin{smallmatrix}1 \\ 1 \end{smallmatrix}\right)\ .$};
\node at (2,1) [fill,circle,inner sep=1.3pt] {};
\node at (2,-1) [fill,circle,inner sep=1.3pt] {};
\end{tikzpicture}
\end{gathered}
\end{equation}
We can thus write the sFD \eqref{reduciblesFD} as $\sigma_2\times 
\sigma_2\times \cdots \times \sigma_2=\sigma_2^{m}$.

Similarly, one has two distinct sFD's with $x=2m-1$ and the two sFD's 
are reducible containing $\sigma_2^{n}$ (for some suitable value of $n$) 
as one of the components and the other component are the following two 
irreducible sFD's that contribute to $c_{1,1}$ and $c_{2,3}$ 
respectively.
\begin{equation}
\begin{gathered}
\begin{tikzpicture}[scale=0.3]
\node at (-3,2) {$\sigma_1\equiv$};
\draw [->] (0,2) -- (5,2) node [below] {};
\node  at (0,2) [fill,circle,inner sep=1.3pt] {};
\node  at (1.5,2) [fill,circle,inner sep=1.3pt] {};
\node  at (3,-1)  {(a)};
\node  at (3,2)  {$\circ$};
\end{tikzpicture}
\qquad\qquad
\begin{tikzpicture}[scale=0.3]
\node at (-4,0) {$\sigma_3\equiv$};
%\node at (10,0) {$\ ,$};
\draw [->] (1,0) -- (-0.7,-1.7) node [left] {};
\draw [->] (1,0) -- (3,-1.4) node [right] {};
\draw [->] (1,0) -- (3,1.4) node [above] {};
\node  at (1,0) [fill,circle,inner sep=1.3pt] {};
\node  at (0.1,-0.9) [fill,circle,inner sep=1.3pt] {};
\node  at (7,0)  {$=\left(\begin{smallmatrix}1&0\\ 1&1\\ 0& 1\end{smallmatrix} \right)\ ,$};
\node  at (3,0)  {$\circ$};
\node  at (1,-1.5)  {$\circ$};
\node at (2,0.7) [fill,circle,inner sep=1.3pt] {};
\node at (2,-0.7) [fill,circle,inner sep=1.3pt] {};
\node  at (2,-4)  {(b)};
\end{tikzpicture}
\end{gathered}\label{sFDzone}
\end{equation}
where we have called the second sFD $\sigma_3$ -- it has two nodes and 
has r.d. $3$. In other words, $c_{2m,2m-1}$ has contributions from two 
sFD's -- one of the form $\sigma_2^{(m-1)}\times \sigma_1$ and the other 
of the form $\sigma_2^{(m-2)}\times \sigma_3$. Studying the symmetries 
of these two sFD's with r.d. $(2m-1)$, one obtains
$$
c_{m,2m-1}=\frac{(2m-1)!}{(2m-2)!!} + \frac{(2m-1)!}{2 (2m-4)!!} =m \times \frac{(2m-1)!}{(2m-2)!!}\ .
$$
Clearly, such a diagrammatic method will enable one to write further 
formulae (we will provide a few more in an appendix) for $c_{m,x}$. 
However, it can get tricky to find all possible diagrams. Keeping this 
in mind, we make the following definition.
\begin{mydef} 
The density, $\rho$, of a sFD with $m$ nodes and r.d. $x$ is $\rho\equiv 
m/x$.
\end{mydef}
The density of a sFD is always greater than or equal to $\tfrac12$ since 
$c_{m,x}=0$ when $x>2m$.
\begin{prop}\label{Dconjecture}
When its density is in the range $(\tfrac12,\tfrac23)$, an sFD  with $m$ nodes and r.d. $x$ is necessarily reducible and one of its components is the sFD, $(\sigma_2)^{n}$, for some $n\geq n_{\textrm{min}}\equiv 2x-3m$.
%as $n_{\textrm{min}}>0$. 
\end{prop} 
The proof follows from Proposition  \ref{Dirreds} that we prove later. When $\rho<2/3$, the 
proposition implies it is impossible to construct an sFD that does not 
contain $\sigma_2$ as a component. The first new sFD, $\sigma_3$, 
appears at $\rho=\tfrac23$. The minimum value of $n$ is fixed by the 
condition that the density of the sFD goes past or equals $\tfrac23$ 
after deleting the nodes that appear in $(\sigma_2)^{n}$ i.e., it is 
smallest value of $n$ such that
$$
\frac{m-n}{x-2n}\geq \frac23 \implies n\geq 2x-3m\ .
$$

\subsection{The third transform}

Proposition \ref{Dconjecture} suggests that in counting the skew FD's 
that contribute to $c_{m,x}$, we can remove components isomorphic to 
$\sigma_2$ in reducible skew FD's and only count skew FD's that do not 
contain any $\sigma_2$ components.  This motivates the next transform 
where we introduce a new triangle $D=(d_{m,x})$.
\begin{equation}
\boxed{
c_{m,x} = \sum_{y=y_{\textrm{min}}}^{m} \frac{x!}{(2y)!! (x-2y)!}\  d_{m-y,x-2y}
}\ , 
\label{Dmatrixdef}  
\end{equation}
with $d_{0,0}=1$, $d_{m,0}=d_{0,m}=0$ for $m>0$ and 
$y_{\textrm{min}}=2x-3m$. The pre-factor in the transform is determined 
by the order of the symmetry of $\sigma_2^{y}$ which is $2^y y!=(2y)!!$.
\begin{enumerate}
\item $d_{m,x}$ counts the number of skew FD's with $m$ nodes and r.d. 
$x$ not containing $\sigma_2$ as its components. Thus it is positive 
definite.
\item Proposition \ref{Dconjecture} implies that $d_{m,x}=0$ when 
$m/x>2/3$. This is stronger than the condition $m/x>1/2$ implied by the 
property of the C-matrix.
\item It is useful to rewrite the transform as follows:
\begin{equation}
c_{m,2m-z} = \sum_{y=\lceil z/2\rceil}^{2z} \frac{(2m-z)!}{(2m-2y)!! (2y-z)!}\   d_{y,2y-z}
\ .
\end{equation}
 In this form, one sees that completely determining row $z$ of the 
D-matrix leads to a nice compact formula for $c_{m,2m-z}$. The 
$D$-matrix clearly contains fewer terms than the $C$-matrix since 
$d_{m,x}=0$ when $\rho<2/3$.

\item To illustrate the transform, consider $c_{m,2m-1}$ which we have 
already computed. One sees that
\begin{align}\label{shiftcformula}
c_{2m,2m-1} &= \sum_{y=1}^{2} \frac{(2m-1)!}{(2m-2y)!! (2y-1)!}\ d_{y,2y-1} \nonumber \\
&= \frac{(2m-1)!}{(2m-2)!!}\  d_{1,1} +  \frac{(2m-1)!}{3!(2m-4)!!} \ d_{2,3} \ .
\end{align}
It is easy to see that $d_{1,1}=1$ as there is precisely one sFD ((a) in 
Eq. \eqref{sFDzone}) and $d_{2,3}=3$ as there are three inequivalent 
diagrams under the action of $S_3$ on the sFD, $\sigma_3$.
\item When $\rho=2/3$, there is only one sFD, $\sigma_3^{m}$, that contributes to 
$d_{2m,3m}$. This implies that
\begin{equation}
d_{2m,3m}= \frac{3m!}{m!\ 2^m}\ , \quad m=1,2,3,\ldots
\end{equation}
\end{enumerate}
$$
D=\begin{pmatrix}
1\\
0 &1  \\
 0& 1 & 3 & 3  \\
0& 1 & 7 & 17 & 28  \\
 0& 1 & 11 & 58 & 156 & 295 & 90  \\
 0& 1 & 18 & 135 & 640 & 1913 & 3786 & 2310  
\end{pmatrix}
$$

\subsection{The final transform}

The main advantage of the $D$-matrix is that it contains fewer terms 
than the $C$-matrix. Using it, we have arrived at formulae for 
$c_{m,2m-z}$ for $z=2,3,4,5$ analogous to the one in Eq. 
\eqref{shiftcformula} that can be obtained, in principle, from the table 
that gives the $D$-matrix. Can we do better? We saw that as the density 
increased from $1/2$ to $2/3$, only one irreducible diagram appears. At 
$\rho=\tfrac34$, two new sFD's appear. They are
\begin{equation}
\sigma_{4a}= \left(\nodee1100\nodee0110\nodee0011\right)\quad,\quad
\sigma_{4b}=\left(\nodee1100\nodee1010\nodee1001\right)\ .
\end{equation}
In fact, one can define another transform that removes reducible 
components of type $\sigma_3$ from sFD;s that contribute to the D-matrix 
for $\rho\in (2/3,3/4)$. The next proposition will enable to do this and 
a lot more by removing a whole family of reducible components that 
necessarily appear when in sFD's with $\rho<1$.

\begin{mydef}
Let $\mathcal{D}\equiv \cup_r \mathcal{D}_r$,  where $\mathcal{D}_r$ denotes the set of strict Ferrers diagrams of dimension 
$r$ consisting only of nodes of the form $(1,1,0,\ldots,0)^T$ or its 
$S_r$ images in addition to the nodes present in $\mu_r$.
\end{mydef}
We say, somewhat loosely, that a strict skew FD, $\sigma$ of r.d. $x$ is 
in $\mathcal{D}$ if the FD $\mu_x\cup \sigma\in \mathcal{D}$. One can 
show that $\sigma_{2}$, $\sigma_3$ and $\sigma_{4a/b}$ are the only 
irreducible strict skew Ferrers diagrams at dimensions $2$, $3$ and $4$ 
respectively that appear in $\mathcal{D}$.

 Let $e_{m,r}$ denote the number of Ferrers diagrams in $\mathcal{D}$ 
obtained by adding $m$ nodes to $\mu_r$. It is easy to see that 
$e_{m,x}=\dbinom{\binom{x}2}m$ as there are $\binom{x}2$ possible nodes 
from which we need to choose $m$ nodes. We define a new transform that 
removes reducible components that are in $\mathcal{D}$. 
\begin{align}
a_{m+r+1,r} &= \sum_{x=1}^r \sum_{p=0}^m \binom{r}{x} e_{m-p,r-x}  \ f_{p+x+1,\ x}\nonumber \\
&= \sum_{x=1}^r \sum_{p=0}^m \binom{r}{x} \binom{\tbinom{r-x}2}{m-p} \ f_{p+x+1,\ x}\ , 
\label{Fdef}
\end{align}
where in the second line we use the explicit formula for $e_{m,x}$ and $f_{1,0}\equiv1$, $f_{n,0}=f_{1,n-1}=0$ for $n>1$. In the first line, a typical term in the summation on the right hand side consists of reducible strict FD's with the  component  in $\mathcal{D}$ having i.d. $r-x$ and $(m-p)$ nodes added to $\mu_{r-x}$ and the other component consisting of an strict FD with no reducible component in $\mathcal{D}$,  i.d. and r.d. $x$ and $p$ nodes added to $\mu_x$ -- their number is counted by $f_{p+x+1,x}$. The binomial factor $\binom{r}{x}$ is the number of ways one can choose $x$ dimensions occupied by the FD's contributing to $f_{p+x+1,x}$.
The 
above formula defines a new triangle $F=(f_{n,r})$.  The entry $f_{r+m+1,r}$ is the
\textit{the number of strict FD's of i.d. $r$ obtained by adding $m$ nodes 
to $\mu_r$ and does \textit{not} contain any reducible components that 
are in $\mathcal{D}$}. Such an FD must necessarily have r.d. also equal to $r$, else it will necessarily have a reducible component isomorphic to $\mu_{r-x}$ if its r.d. is $x$.

It is easy to see that $f_{r+1,r}=0$. The only contribution to 
$a_{r+1,r}$ is the unique FD $\mu_r$ which is $\mathcal{D}$. Similarly, 
$f_{r+2,r}=0$ when $r>1$ as the only contribution to $a_{r+2,r}$ is of the 
form $\sigma_1\times \sigma_2^{r-1}$. One also has $f_{3,1}=1$ with $\sigma_1$ being the unique
FD contributing to it. The next proposition shows the advantage of defining the $F$-matrix.

\begin{prop}\label{Fprop}
$f_{m+r+1,r}=0$ when $r>m$.
\end{prop}
\textbf{Proof:} Let $\lambda$ be an FD of i.d. $r$ with $m+r+1$ nodes 
that contributes to $f_{m+r+1,r}$. Consider the skew FD, 
$\lambda\setminus\mu_r$ -- it has $m$ nodes. It must be a strict skew FD 
else it has a irreducible component isomorphic to $\mu_x$ for some 
$x<r$. Thus, the proposition implies that there are no strict skew FD's 
with density $\rho=m/r<1$.

We can also assume that the skew FD is irreducible -- if it is 
reducible, it must necessarily have at least one irreducible component 
with density $<1$ and we can focus on (proving the non-existence) such 
irreducible components. Our goal is thus reduced to proving that there 
are no irreducible strict skew FD's with density $<1$.

\begin{mydef}
Let us call the nodes obtained by all permutations of the coordinates of 
the node $(1,1,0,\ldots,0)^T$ as nodes of type $1$. Similarly, call the 
nodes obtained by permuting coordinates of $(2,0,\ldots,0)^T$ as type 
$2$. Nodes of type $3$ are nodes that are not of type $1$ or $2$.
\end{mydef}
Examples of type 3 nodes include $(1,1,1,0,\ldots,0)^T$ and 
$(3,0,\ldots,0)$. Such nodes cannot be added to the FD $\mu_r$ without 
including supporting nodes of type $1$ and $2$. The addition of nodes of 
type $3$, when possible, never increase the r.d. of an FD thus 
increasing the density. Thus, given a FD $\lambda$ (of i.d. $r$ and r.d. 
$r$) containing type $3$ nodes, we can form a new FD $\lambda'$ with the 
same r.d. but lower density. Further, if $\lambda\setminus\mu_r$ is 
irreducible, $\lambda'\setminus\mu_r$ is also irreducible. The skew FD 
$\lambda'\setminus\mu_r$ thus consists of nodes of type $1$ and type 
$2$. If it consists of only nodes of type $1$, then 
$\lambda'\in\mathcal{D}$. Thus, we need to only consider irreducible 
strict skew FD's containing at least one node of type $2$.

For the rest of the discussion, let $\lambda'$ be an FD such that 
$\lambda'\setminus\mu_r$ is an irreducible strict skew FD containing 
only nodes of type $1$ and at least one node of type $2$. It is easy to 
see that removing of node of type $2$ does not affect  the 
irreducibility of the skew FD. Further, it does not reduce the r.d. as 
the only way a type $2$ node can reduce the r.d. of a skew FD is when it 
appears as a part of a reducible component isomorphic to $\sigma_1$. 
Thus, we can delete all type $2$ nodes to obtain a new FD $\lambda''$ 
that is irreducible and contains only type $1$ nodes. Again, it is easy 
to see that $\rho(\lambda'')\leq \rho(\lambda')$. Further 
$\lambda''\in\mathcal{D}$. Thus, one has the sequence
\begin{equation}
\rho(\lambda'')\leq \rho(\lambda')\leq \rho(\lambda)\ .
\end{equation}
Let $\lambda'\setminus \mu_r$ have $(r-1)$-nodes so that its density is 
just below one and contain $z$ nodes of type $2$. Then, 
$\lambda''\setminus \mu_r$ will have $(r-1-z)$ nodes and be irreducible. 
The next proposition shows that such a $\lambda''$ does not exist.  
Hence, there exists \textit{no} FD $\lambda''$ and hence no FD $\lambda'$ 
with density $<1$.\mbox{}\hfill$\Box$

\begin{prop}\label{Dirreds}
The only strict FD's in $\mathcal{D}$ of i.d. $r$ such that the skew FD 
$\lambda\setminus\mu_r$ is strict and irreducible with density less than 
$1$ necessarily have $\rho=\tfrac{r-1}r$.
\end{prop}
\textbf{Proof:}  Let us assume 
that $\lambda\setminus \mu_r$ has $(r-2)$ nodes and is irreducible. Let 
us try to construct such a strict skew  FD and we will see that there are not enough 
nodes. Start by putting the first type $1$ node in the $x^1x^2$ plane. The 
irreducibility condition implies that the second node must be either in 
the $x^1x^\alpha$ or $x^2x^\alpha$ plane where $\alpha$ is not $1$ or 
$2$. The key point is that the additional node must contain one of the 
used up coordinates, $x^1$ or $x^2$ in this case and a new coordinate so 
that irreducibility is maintained. Clearly, such a process needs $(r-1)$ 
nodes to get an irreducible skew FD $\lambda\setminus\mu_r$ with r.d. 
$r$. This is impossible. Hence, there exists \textit{no} irreducible 
skew FD $\lambda$ with density $\tfrac{r-2}r$.  it is easy to extend the 
argument to exclude even lower densities. Thus, the only possibility that is not
ruled out is to have strict skew FD's with $(r-1)$ nodes with r.d. $r$ -- these have
density $\tfrac{r-1}r$. 
\\ \mbox{}\hfill $\Box$

\noindent \textbf{Remark:} $\sigma_2$, $\sigma_3$ and $\sigma_{4a/b}$ are the only irreducible
strict skew FD's with r.d. $2,3,4$ respectively.

\subsubsection{Properties of the $F$-matrix}
\begin{enumerate}
\item The most important property is the one implied by Proposition \ref{Fprop} which says that the $F$ matrix is lower triangular with $f_{n,r}=0$ when $r<[(n-1)/2]$. For fixed value of $n$, the $F$-matrix has far fewer terms (roughly half) than the corresponding row in the $A$-matrix. 
We have determined the first $25$ rows of the $F$-matrix (see Table \ref{Ftable}).
\item It turns out that there are other transforms that also lead to matrices with fewer entries like the $F$-matrix. See for instance, the box transform that we consider in the appendix. However, their relationship to $A$ is not as simple as Eq. \eqref{Fdef}. The simplicity of Eq. \eqref{Fdef} is what picks out the $F$-matrix as special.
\item We can also use this idea to refine the counting problem associated with the $C$-matrix. 
Let $C^{\mathcal{D}}=(c^{\mathcal{D}}_{m,x})$ denote the contribution to the $C$-matrix that arise from FD's that are in $\mathcal{D}$. Since the set $\mathcal{D}_r$ is invariant under $S_r$, it is easy to see that $C^{\mathcal{D}}$ is given by the transform
\begin{equation}
 \binom{\tbinom{x}2}{m}= \sum_{x=0}^{2m} \binom{r}{x} c^{\mathcal{D}}_{m,x}\ . 
 \end{equation}
Then, we can define $\widetilde{C}=(\widetilde{c}_{m,x})$ by removing contributions that arise from reducible parts that are isomorphic to contributions to $\mathcal{C}^{\mathcal{D}}$. Then, one has
\begin{equation}
c_{m,x} = \widetilde{c}_{m,x} +c^{{\mathcal D}}_{m,x}+\sum_{y=1}^{x-1} \sum_{p=1}^{m-1}\binom{x}{y}\  c^{{\mathcal D}}_{m-p,x-y}\  \widetilde{c}_{p,y}\  .
\end{equation}
Given a strict skew FD that contributes  to $ \widetilde{c}_{m,x} $, it is easy to see that there is a unique FD obtained by adding nodes in $\mu_x$ to the skew FD. Further,  this FD must contribute to the entry $f_{m+x+1,x}$ in $F$. Since the converse also holds i.e, given a strict  FD of i.d. $x$ that contributes to the $F$-matrix, the skew FD obtained by deleting nodes in $\mu_x$ gives a skew FD that contributes to $\widetilde{C}$. Thus, one has 
\begin{equation}
\widetilde{c}_{m,x} =f_{m+x+1,x}\ .
\end{equation}
\end{enumerate}

We observe numerically that $f_{2m+1,m}=(m+1)^{m-2}$ for $m=0,1,2,\ldots,12$. We will show  that it holds for all $m$. These numbers appear in the  sequence numbered A000272 in the OEIS\cite{oeis}. The next proposition presents a further refinement. We need a few definitions which we briefly state. A graph, consisting of vertices and undirected edges, with no cycles is called an acyclic graph or a forest. A forest may consist of disconnected components and is called a tree if it has only one connected component. A rooted tree is one with a marked/special vertex (called the root) while a rooted forest is one in which every component is rooted. A spanning forest is any subgraph that is both a forest (contains no cycles) and spanning (includes every vertex)\cite{wiki:span tree,Flajolet:2009}.
\begin{prop}\label{Cayley}
Let $\alpha$ be the number of nodes of type $2$ contained in an  FD that contributes to $f_{2m+1,m}$.  Let $f_{2m+1,m}(\alpha)$ denote the total number of such Ferrers diagrams. Then, $f_{2m+1,m}(\alpha)$ is the number of spanning rooted forests on $m$ vertices and $\alpha$ components. It follows from a result due to Cayley on the numbers of spanning rooted forests that\cite{Pak:2009}
\begin{equation}
f_{2m+1,m}(\alpha)=\binom{m-1}{\alpha-1}\ m^{m-\alpha}\ .
\end{equation}
\end{prop}
\textbf{Proof:}
We will provide  a bijective map relating FD's that contribute to  $f_{2m+1,m}(\alpha)$  to spanning
rooted forests on $m$ vertices and $\alpha$ components. There is a natural action of $S_m$ on both sides -- on the FD side, it corresponds to permuting the $m$ coordinates and on the rooted forest side, it corresponds to relabeling the $m$ nodes. We  identify these two groups. 

Given a skew FD that contributes to $f_{2m+1,m}(\alpha)$, we can construct a graph with $m$ vertices labeled from $(1,\ldots,m)$ as follows. The type $2$ nodes become root vertices carrying the label of the non-vanishing coordinate. Thus if a type $2$ node has non-vanishing $j$-th coordinate, assign it the label $j$. Add $(m-\alpha)$ vertices and label them with the unused labels. Every type $1$ vertex has two non-vanishing coordinates, say the $j$-th and $k$-th coordinates. Assign an edge that connects vertex $j$ to vertex $k$. Repeat for all type $1$ nodes. In this process, there are as many components as there are type $2$ nodes. 
Thus the graph is a spanning rooted forest on $m$ vertices and $\alpha$ components.  The following example illustrates the map for $m=4$ and $\alpha=1$. The root vertex is shown by a filled circle.
\begin{center}
%\boxed{
%\begin{tikzpicture}[scale=0.4]
%\node at (-6,0.5) {$\left(\begin{smallmatrix}
%2 & 1 & 0 & 0 \\
%0 & 1 & 1 & 0 \\
%0 & 0 & 1 & 1 \\
%0 & 0 & 0 & 1 
%\end{smallmatrix}\right)$};
%\node  at (0,2) [fill,circle,inner sep=1.5pt] {};
%\node at (0,2) [right] {{\tiny 1}};
%\draw[-] (0,2)--(0,1.2);
%\node at (0,1) {$\circ$};
%\node at (0,1) [right] {{\tiny 2}};
%\draw[-] (0,0.8)--(0,0.2);
%\node at (0,0) {$\circ$};
%\node at (0,0) [right] {{\tiny 3}};
%\draw[-] (0,-0.2)--(0,-0.8);
%\node at (0,-1) {$\circ$};
%\node at (0,-1) [right] {{\tiny 4}};
%\end{tikzpicture}
%}
%
%\boxed{
\begin{tikzpicture}[scale=0.5]
\node at (-8,1) {$\left(\begin{smallmatrix}
2 & 1 & 1 & 0 \\
0 & 1 & 0 & 0 \\
0 & 0 & 1 & 1 \\
0 & 0 & 0 & 1 
\end{smallmatrix}\right)$};
\node at (-4,1) {$\longleftrightarrow$};
\node  at (0,2) [fill,circle,inner sep=1.5pt] {};
\node at (0,2) [right] {{\tiny 1}};
\draw[-] (0,2)--(-1,1.2);
\node at (-1,1) {$\circ$};
\node at (-1,1) [right] {{\tiny 2}};
\draw[-] (0.2,1.9)--(0.9,1.2);
\node at (1,1) {$\circ$};
\node at (1,1) [right] {{\tiny 3}};
\draw[-] (1.1,0.9)--(1.9,0.1);
\node at (2,0) {$\circ$};
\node at (2,0) [right] {{\tiny 4}};
\end{tikzpicture}
%}
\end{center}

To prove the converse statement, given a spanning rooted forest with $m$ vertices and $\alpha$ components, we need to construct an FD that contributes to $f_{2m+1,m}(\alpha)$. This is easy to do . Pick the root vertices  and assign them to type $2$-nodes whose non-vanishing coordinate decided by the label of the vertex. Next assign to all edges a type $1$ node that has non-vanishing coordinates at precisely the locations decided by the labels of the vertices it connects.  We thus recover the skew FD.
\mbox{}\hfill $\Box$\\

\noindent\textbf{An example:} We know that $f_{5,2}=3$. The three skew FD's are: 
$$
\sigma_1^2 =\left(\begin{smallmatrix} 2 & 0 \\ 0 & 2 \end{smallmatrix}\right)\quad;\quad
\left(\begin{smallmatrix} 2 & 1 \\ 0 & 1 \end{smallmatrix}\right)\quad;\quad
\left(\begin{smallmatrix} 1 & 0 \\ 1 & 2 \end{smallmatrix}\right)\quad,
$$
where the nodes are listed by the ordering: $(a_1,a_2)> (b_1,b_2)$ if $a_1>b_1$ or $a_1=b_1$ and $a_2>b_2$.
Note that there are two equivalence classes of skew FD's. Under $S_2$ action as the second and third skew FD's get mapped to each other.
\begin{center}
\begin{tikzpicture}[scale=0.4]
\node  at (0,0) [fill,circle,inner sep=1.5pt] {};
\node at (0,0) [above] {{\tiny 1}};
\node  at (2,0) [fill,circle,inner sep=1.5pt] {};
\node at (2,0) [above] {{\tiny 2}};
\node at (4,0) {;};
\draw[-](6,0) -- (8,0);
\node at (6,0) [above] {{\tiny 1}};
\node  at (6,0) [fill,circle,inner sep=1.5pt] {};
\node  at (8.2,0)  {$\circ$};
\node at (8.2,0) [above] {{\tiny 2}};
\node at (10,0) {;};
\draw[-](12,0) -- (14,0);
\node at (11.8,0) [above] {{\tiny 1}};
\node at (14,0) [above] {{\tiny 2}};
\node  at (14.1,0) [fill,circle,inner sep=1.5pt] {};
\node  at (11.8,0)  {$\circ$};
\node at (17,0) {.};
\end{tikzpicture}
\end{center}
\textbf{Remark:} Given a skew FD, it is possible to uniquely label the nodes by ordering them by a choice of ordering as illustrated above.

\section{Other triangles}

\subsection{New triangles}

So far, we have considered transforms that lead to new triangles 
($A/C/D/F$) all of which have positive definite entries since we they 
all count numbers of skew Ferrers diagrams. We will now provide two 
other transforms that are partly conjectural and lead to triangles that 
are not positive definite -- we denote the entries with Greek letters to 
remind us of this. We begin by expanding the entries in the $A$-matrix 
as follows. Let
\begin{equation}
 a_{m+r+1,r} = \sum_{z=0}^{2m} \alpha_{m,z}\  \frac{r^{2m-z}}{2m!!}\ ,
\end{equation}
with $\alpha_{m,0}=1$ for $m\geq0$ and $\alpha_{m,2m}=0$ for $m>0$. The 
above transform provides the entries for another triangular matrix, 
$\alpha_{m,z}$, that we call the $\alpha$-triangle by setting 
$\alpha_{m,z}=0$ for $z>2m$. One can explicitly relate the 
$\alpha_{m,z}$ to the entries in the $C$-matrix using Stirling numbers 
of the first kind. Thus the above formula is not conjectural. However 
the following is conjectural:
\begin{conj}
The entries of the $\alpha$-triangle, i.e., $\alpha_{m,z}$, are all 
integers.
\end{conj}
This is true for the first ten rows and appears to hold for the first eleven rows which have been determined using conjectures.

The second conjecture introduces a new triangle, that we call the 
$\beta$-triangle, and its associated transform. It has been determined 
experimentally and verified to hold to the extent possible.
\begin{conj} The $\alpha$-matrix admits the following decomposition.
\begin{equation}
\alpha_{m,z}= \sum_{y=0}^{\lfloor z/2\rfloor} \binom{m}{z-y}\ \beta_{z,y}\   ,
\end{equation}
with $\beta_{0,0}=1$ and $\beta_{2y,y}=0$ for all $y>0$.  
\end{conj}
By setting $\beta_{z,y}=0$ for $y>\lfloor z/2\rfloor$, this becomes the 
binomial transform
\begin{equation}
\alpha_{m,z}= \sum_{y=0}^{m} \binom{m}{z-y}\ \beta_{z,y}\   .
\end{equation}
The inverse transform is
\begin{equation}
\beta_{z,y}=\sum_{m=0}^{z-y} (-1)^{m+z-y} \binom{z-y}{m}\ \alpha_{m,z}\ .
\end{equation}
We now state a conjecture of Meeussen that fixes one of the coefficients. 
\begin{conj}[Meeussen]
$$
\beta_{n,0}=H_n(\tfrac12)\ ,
$$
where $H_n(x)$ is the $n$-th Hermite polynomial.
\end{conj}
Recall that the $\alpha$-matrix has $2m$ non-zero entries in the $m$-th 
row. The $\beta$-matrix has fewer terms, roughly half  the entries  in
the $\alpha$-matrix. We were able to determine eleven rows of the $\alpha$ and $C$-matrices using the $\beta$-matrix of which $10$ were verified through other means. This was our main motivation
in searching for and find the combinatorial problem that eventually lead to the $F$-matrix.

\subsection{The B-triangle}

We now construct another lower triangular matrix $B=(b_{n,r})$ with 
$n=1,2,\ldots$ and $r=0,1,2,\ldots$ and $b_{n,0}=1$.
\begin{equation}
p_d(n) = \sum_{r=0}^{n-1} \binom{d}{r} \ b_{n,r}=1+ \sum_{r=1}^{n-1} \binom{d}{r} \ b_{n,r}\ . \label{bdef}
\end{equation}
The  matrix $B$ appear in the OEIS as sequence number A096806.
Using Pascal's identity
\begin{equation}
\binom{d+1}r = \binom{d}r + \binom{d}{r-1}\ ,
\end{equation}
we can relate the matrix $B$ to $A$. Thus, one has the relation
\begin{equation}
b_{n,r}=a_{n,r}+a_{n,r+1}\ .
\end{equation}
One can easily show that $b_{n,n-1}=1$ using the above formula and known 
properties of the matrix $A$. The first six rows of $B$ have been 
determined explicitly, for instance, in Andrews' book on 
Partitions\cite{AndrewsPartitions}. It is easy to check that the above 
relation holds for all six rows.

\subsection{Hanna's triangle}

\begin{conj}[Hanna]
There exists a lower-triangular matrix $T=(\tau_{ij})$ (with 
$i,j=0,1,2,\ldots$) with integral entries and ones on its diagonal such 
that
$$
p_d(n) =\sum_{j=0}^n (T^d)_{n,j}\ .
$$
\end{conj}
In other words, the sum of the the $n$-th row of the $d$-th power of $T$ 
give the $d$-dimensional partition of $n$. This matrix appears in the 
OEIS as sequence A096651.  Since $p_d(0)=1$, we can set $\tau_{0,0}=1$ 
and $\tau_{j,0}=0$ for $j>0$. For the rest of the discussion, we will 
consider $n>0$ and can delete the zeroth row and column of the 
$T$-matrix as they no longer play a role. We shall however use the same 
symbol $T$ to denote the modified matrix as it is easy to reconstruct 
the original $T$ matrix by adding back the zeroth row and column. We 
shall prove the existence as well as the integrality of the matrix $T$ 
by constructing an explicit map that relates $T$ to the matrix $B$ (and 
hence $A$) that we considered in the previous section.

\noindent\textbf{Proof:} For $n\geq 1$, the Hanna conjecture can written as
\begin{equation}
p_d(n) = \sum_{j=1}^{n} (T^d)_{n,j}=   \sum_{x_1\cdots x_{d}} \tau_{n,x_1}\tau_{x_1,x_2}\cdots \tau_{x_{d-1,}x_d}\ ,
\end{equation}
where $ n\geq x_1 \geq x_2\geq \cdots \geq x_d\geq 1 $. It obviously 
holds for $n=1$ since $\tau_{11}=1$. Using the fact that $ T $ has ones 
in its diagonal, we can simplify the above expression to
\begin{equation}
p_d(n) = 1+\sum_{r=1}^{n-1} \binom{d}{r}  \sum_{x_1\cdots x_{r}} \tau_{n,x_1}\tau_{x_1,x_2}\cdots \tau_{x_{r-1},x_r}\ .
\end{equation}
with sum now running over all sequences of $r$ positive non-zero 
integers $ (x_1,\cdots, x_{r}) $ such that $ x_0\equiv n> x_1 > x_2> 
\cdots > x_r \geq 1$. The combinatorial factor expresses the number of 
ways in which diagonal elements are chosen. Comparing the above equation 
with Eq. \eqref{bdef} implies the (potential) identity for $n>1$ and 
$r\geq1$.
\begin{equation}
\sum_{x_1\cdots x_{r}} \tau_{n,x_1}\tau_{x_1,x_2}\cdots \tau_{x_{r-1},x_r}= b_{n,r}\ ,
\end{equation}
with $n> x_1 > x_2> \cdots > x_r \geq 1$. Let us assume that this 
relation holds for $n<m)$ (for some $m>1$) and that we have determined 
$(m-1)$ rows of $T$. Then, we can rewrite the above equation as
\begin{equation}
\sum_{1\leq x<m} \tau_{m,x} b_{x,(r-1)} = b_{m,r}\  \textrm{ for } m>r\geq1 \  . \label{tsolve}
\end{equation}
The above $(m-1)$ equations are linear equations in $(m-1)$ unknowns: \\ 
($\tau_{m,1},\ldots, \tau_{m,m-1})$ -- these are the undetermined 
entries in the $m$-th row of $T$. Hence, they have a solution if the 
matrix (constructed using $b_{x,(r-1)}$) is invertible. The matrix is 
upper triangular with ones in its diagonal. Hence it is has determinant 
one and hence is invertible. This enables us to recursively determine 
all the entries in the matrix $T$. This proves the \textit{existence} of 
$T$.

We shall inductively prove the integrality of the matrix $T$ using more 
explicit details of Eq. \eqref{tsolve}. We begin with the equation for 
$r=m-1$ and it gives
\begin{equation}
\tau_{m,m-1} b_{m-1,m-2} = b_{m,(m-1)} \implies\boxed{\tau_{m,m-1} = 1 }\ ,
\end{equation}
where we have used $b_{m,m-1}=1$ for $m\geq1$. Next consider, $r=m-2$. 
This equation gives $\tau_{m,m-2} +\tau_{m,m-1} b_{m-1,m-3} = b_{m,m-2}$ 
which gives
\begin{equation}
\tau_{m,m-2} = b_{m,m-2}-\tau_{m,m-1} b_{m-1,m-3} \ ,
\end{equation}
where we have used the fact that $\tau_{m,m-1}$ has been solved for and 
shown to be integral in the previous step. Note that this implies that 
$\tau_{m,m-2}$ is integral.  Proceeding in this manner from $r=(m-1)$ to 
$r=1$, we thus determine all the unknowns. A typical equation will take 
the form (reflecting the triangular nature of the equations)
\begin{equation}
\boxed{\tau_{m,m-r}= b_{m,m-r} -\sum_{x=m-r+1}^{m-1} \tau_{m,x} b_{x,m-r} }\ ,
\end{equation}
for $r=1,2,\dots, (m-1)$. We assume that $\tau_{m,m-r'}$ is integral for 
all $r'<r$. Thus the right hand side is integral as it only contains 
integral terms.  Hence $\tau_{m,m-r}$ is integral. This concludes the 
proof of integrality of the matrix $T$. \\ \mbox{}\hfill $\Box$

\noindent We now state an unproven conjecture of Hanna and Meeussen.
\begin{conj}[Hanna-Meeussen]
$m!\ \tau_{m+r+1,m}$ is a polynomial of degree $m$ in $r$ with integral 
polynomial coefficients.
\end{conj}
It is easy to show that $\tau_{m+r+1,m}$ is a polynomial of degree 
$(2m-1)$ in $r$ using the properties of the $A$-matrix. However, the 
above conjecture is stronger and seems to consistent with known data for 
$m=0,1,\ldots,11$.

\section{Practical Considerations}

This section provides details on the exact enumeration of 
higher-dimensional partitions as well as the triangles defined in this paper. With access to high-performance computing getting 
easier in recent times, this is indeed an additional computational 
aspect that can and must be added to the theoretical discussion of the 
previous section. We will first discuss the algorithms that we used and 
then discuss exact enumerations as we carried out.

\subsection{Algorithms for higher-dimensional partitions}

There are two algorithms in the literature for computing 
higher-dimensional partitions. The first one is due to Bratley and McKay 
(the BM algorithm)\cite{Bratley:1967a} and the second one is due to 
Knuth\cite{Knuth:1970} -- both are more than 40 years old reflecting the 
lack of progress in this area. Both are highly recursive and provide 
distinct ways of exactly enumerating higher dimensional partitions.

\subsubsection*{The BM algorithm}

The partitions in any fixed dimension, say $d$, form a tree which we 
call the \textit{partition} tree in $(d+1)$-dimensions\footnote{Recall 
that the Ferrers diagram for a $d$-dimensional partition is a set of 
points in $d+1$ dimensions.} and denote by the symbol 
$\mathcal{T}_{d+1}$. Every node of the tree is the Ferrers diagram 
associated with a partition. The unique Ferrers diagram containing one 
point is the root node of the tree. New partitions can be formed by 
adding or deleting a point from the Ferrers diagram\footnote{To avoid 
confusion, in this section alone, we shall refer to nodes of a partition 
as points in the Ferrers diagram. This is to avoid confusion with the 
node of the tree.}. Add a link to partitions connected this way. The 
depth of the tree is the number of points in the partition.

The BM algorithm recursively traverses the tree up to some fixed depth, 
say $n$, such that each node is visited precisely once. The heart of the 
algorithm is the routine called $\textit{part}$ that takes three 
arguments and is recursively called in the algorithm. Every time a node 
is visited, the partition is stored in an array called \textit{current} 
and presented to user. If one is interested in only counting the number 
of partitions of an integer in a given dimension, if the current 
partition has $m$ points, increment a suitable counter, call it 
$p_d(m)$, by one. At the end of the program, the counter thus contains 
the number of partitions of all integers less than or equal to the depth 
of the traversed tree.

\subsubsection*{The Knuth algorithm}

Let $S_m=\mathbb{N}^m$ denote the set of points in the totally positive 
orthant in a hyper cubic lattice. Let $d_m(k)$ denote the number of 
topological sequences with index $k$ (see 
\cite{Knuth:1970,Balakrishnan:2011bm} for definitions). Then a theorem 
due to Knuth\cite{Knuth:1970} relates the numbers of topological 
sequences to numbers of partitions. To be precise, one has
\begin{equation}
p_m(n)=\sum_{k=0}^n d_m(k) p_1(n-k)\ .
\end{equation}
Since one-dimensional partitions are easily enumerated from the 
generating function, it is simple to generate $p_m(n)$ given $d_m(k)$ 
for all $k\leq n$. Knuth provided an algorithm to generate and count all 
topological sequences -- he illustrated this method by generating 
numbers for the numbers of solid partitions for integers $\leq 28$. 
Recently, a parallelized version of this algorithm was used by the 
author and other collaborators to enumerate solid partitions of integers 
$\leq 68$\cite{Balakrishnan:2011bm}. \\

\noindent\textbf{Remark:} An important aspect of the BM algorithm is 
that its memory usage is of the order of $nd$ bytes, where $d$ is the 
dimension and $n$ is the maximum depth. This is vastly superior to the 
Knuth algorithm, where a similar problem needs memory of the order of 
$n^{d-1}$ bytes. However, when memory isn't an issue, our implementation 
of the Knuth algorithm typically takes less time than our implementation 
of the Bratley-McKay algorithm.

\subsubsection*{The modified BM algorithm}

We begin with the observation that a suitably chosen sub-tree of the 
partition tree in $r$-dimensions, $\mathcal{T}_r$ generates all 
partitions that contribute to the $r$-th column of the $A$-matrix i.e., 
$a_{n,r}$. The head node of this sub-tree is the Ferrers diagram $\mu_r$ 
defined in Eq. \eqref{specialFD}. The rest of the tree is generated by 
adding points to $\mu_r$. Let us denote this sub-tree by $\mathcal{V}_r$ 
and the depth of this tree is clearly $m$ where $m=n-r-1$.

The BM algorithm was designed to recursively traverse the partition tree 
visiting each node precisely once. The starting point of the algorithm 
is the root node whose Ferrers diagram consists of one point. Our idea 
is to change the initial configuration in the BM algorithm to the 
Ferrers diagram, $\mu_r$ and then call the recursive routine 
\textit{part} with suitably chosen arguments\footnote{We have determined 
that the correct call is $part\left(r+2,0,\tbinom{r+1}2\right)$. For 
comparison, the BM algorithm begins with the call $part(1,0,1)$. We 
thank Arun K. Jayaraman for implementing the BM algorithm as well as 
working out this modification. }. For this modification to work 
correctly, the program should traverse the sub-tree $\mathcal{V}_r$ 
visiting each node precisely once to the chosen depth. This turned out 
to be easier as we experimentally observed that the sub-tree 
$\mathcal{V}_r$ appeared naturally in the original BM algorithm for low 
values of $r$. We then checked that the modified BM algorithm correctly 
generated entries in the $A$-matrix for $r\leq 10$. However, we have 
\textit{not} rigorously proved that this is indeed the case.

Thus, once we have the modified BM algorithm correctly traversing the 
sub-tree $\mathcal{V}_r$, we can do the following:
\begin{itemize}
\item Count the number of nodes at each depth -- this gives the number 
$a_{m+r+r,r}$.
\item At each node, numerically compute the reduced dimension, $x$ of 
the Ferrers diagram. Then organizing the partitions by depth and r.d., 
we determine $\binom{r}{x} c_{m,x}$. The binomial pre factor is present 
since all $x\leq r$ will appear. This also implies that the algorithm is 
inefficient computationally for obtaining entries in the $C$-matrix.
\end{itemize} 

\subsubsection*{A wish list of algorithms}

As we just mentioned, the current algorithm to enumerate entries in the 
$C$-matrix is computationally inefficient as we generate $\binom{r}{x}$ 
partitions for each distinct contribution to $c_{m,x}$. It is also 
inefficient because we need to compute $x$ for every given partition. 
Can we create a more efficient algorithm? The problem is that we do not 
have an elegant characterization of sFD's with r.d. equal to $x$. This 
is in contrast to what happened with the $A$-matrix. In that case, we 
could show that any FD that has i.d. $r$ necessarily contains the FD 
$\mu_r$. By using it as our initial configuration, we directly avoided 
configurations with smaller intrinsic dimension. For the $C$-matrix, we 
cannot avoid configurations that have smaller r.d. than the one of 
interest.

We do not have any algorithms for the $\alpha$ and $\beta$ matrices as 
well as the $D/F$ matrices. So far these have been computed only 
indirectly after the $A$ and $C$ matrices have been computed. However, Proposition \ref{Cayley} might be a good starting point to coming up with an algorithm that
 directly enumerates entries in the $F$-matrix.

\subsection{Exact enumeration of higher-dimensional partitions}

In order to evaluate higher-dimensional partitions for integers $\leq 
25$ and dimensions $\leq 10$, we chose to use the Knuth algorithm do 
carry out our computations. There was no serious memory issues for 
dimensions $\leq 7$ and the Knuth algorithm worked well.

We needed to modify our computation when for dimensions $8$, $9$ and 
$10$. The reduction in memory was done by counting topological sequences 
that fit into a box of size $b$. Then the memory requirement went down 
from $n^{d-1}$ to $b^{d-1}$. For instance, when $n=20$ and $b=10$ (for 
$d=10$), the memory usage went down by a factor of $2^9$ and enabled us 
to keep our memory requirements in the $4-8$ GB range as constrained by 
the IITM supercluster. However, some configurations are missed out as 
they do not fit into the box. Interestingly, one can show the error due 
to missed configurations is independent of box size when the index lies 
in the range $[b+1,2b]$. This makes it easy to estimate the errors by 
comparing with known results at smaller values of $b$ and then slowly 
increasing the value of $b$. This method was used, for instance, to 
determine the ten-dimensional partitions of $20$ -- this was carried out 
by using a box of size $11$ with errors determined up to $k=b+9$. This 
was one of the more difficult computations as it took a several months 
of computer time to first estimate the errors and then carry out the 
final run in the box. Table \ref{Knuthresults} gives the results 
obtained used the Knuth algorithm for $n\leq 23$ and $d\leq10$ and 
represent more than six months of computer time.

\subsection{Exact enumeration of the $A$ and $C$ triangles}

The modified BM algorithm was used to generate the $A$ and $C$ matrix.  
The first eight rows of the $C$-matrix have been completely determined. Two additional rows
were determined using additional information from the $D$-matrix. We obtain
\begin{align}
c_{m,2m-2} &= \tfrac{(2 m - 2)! }{6(2m-4)!!} (3 m^2 - m - 1) \nonumber \\
c_{m,2m-3} &=\tfrac{(2m-3)!}{6(2m-4)!!}(2 m^4 - 6 m^3 + 3 m^2 + 3 m + 4)\\
c_{m,2m-4} &=\tfrac{(2 m - 4)!}{180 (2 m - 6)!!}(15 m^5- 75 m^4+ 95 m^3 + 21 m^2  + 88 m   + 42) \nonumber\\
c_{m,2m-5} &=  \tfrac{(2 m - 5)!}{ 90 (2 m - 6)!!} (258 - 167 m - 80 m^2 + 111 m^3 - 174 m^4 + 
   116 m^5 - 31 m^6 + 3 m^7) \nonumber
\end{align}
This determines all entries in the $A$-matrix of the form $a_{m+r+1,r}$ 
for $m=0,\ldots 10$ for \textit{all} values of $r$. We have determined 
the remaining entries for $a_{n,r}$ for $n\leq 23$ by using the BM 
algorithm when necessary. The entry $a_{23,11}$ was one of the longest 
runs and took about 880 hours of CPU time. Tables \ref{Amatrix} and 
\ref{Cmatrix} provide our results.

Using 
the $\beta$ and $\alpha$ matrices as well as the Meeussen conjecture, we 
have also determined the $11$-th row of the $C$-matrix. While none of 
these results were used in finally determining the entries in the 
$A$-matrix, there doesn't seem to be an inconsistency. This is only to 
be viewed as evidence for various conjectures.

\subsection{Extracting the elements of the other triangles}

All other triangles were obtained by using known numbers for the $A$ and 
$C$ matrices as we do not have an algorithm to enumerate them. The 
results for the $D$  and the $\beta$-matrix are presented in Tables 
\ref{Dmatrix} and \ref{Betamatrix} respectively.

An improved implementation of the Bratley-McKay algorithm was provided to us recently by Prof. Bratley. This enabled us to enumerate a few more terms -- in particular, we were able to enumerate rows $24$ and $25$ up to and including $a_{25,12}$. This enabled us to completely determine $25$ rows
of the $F$-matrix. This in turn determines all entries in $25$ rows of the $A$-matrix and hence determines partitions of $25$ in any dimension. It also provides a check on the 23 rows of the $A$-matrix which was independently determined. Table \ref{Ftable} provides our results.

\section{Concluding Remarks}

We have shown the existence of several structures that lead to 
simplifications in the exact enumerations of higher-dimensional 
partitions. The combinatorial interpretations that we have provided have 
enabled us to come up with an algorithms to evaluate the $A$ and $C$ 
matrices.   A few lines 
of code in Mathematica/Maple/Maxima/java can be used to store the $A$ 
matrix and compute $p_d(n)$ for $n\leq25$ using the Binomial transform 
in real time. A working implementation of this is provided on the webpage:
\begin{center}
 \texttt{http://www.physics.iitm.ac.in/\symbol{126}suresh/partitions.html}\ .
 \end{center}
We will be adding these numbers to the OEIS as well as 
providing modules for SAGE/Mathematica/Maxima.

It appears difficult to improve on our results which have determined all 
entries for the $n=25$ row of the $A$-matrix. In fact, we have 
determined most of the entries for the $n=26$ entry and hope to add this 
 row in the future. Further additions to the $A$-matrix will require new and efficient 
 algorithms  to directly enumerate either the $C$ or the $F$ matrix.
Another approach would be a na\"ive parallelization of the BM algorithm.  We hope to be able to 
eventually determine partitions of integers less than $30$ in any 
dimension in the future.
\bigskip

\bigskip

\noindent

\textbf{Acknowledgments:} This work arose out my attempt to prove 
Hanna's conjecture. I thank Paul Hanna for drawing my attention to his 
conjecture. Paul Hanna and Wouter Meeussen have shown a passion for 
sequences that is highly infectious and have been a constant source of 
encouragement during my attempt at enumerating the various triangles. I 
am grateful to Arun K. Jayaraman, then an undergraduate student, who 
implemented the Bratley-McKay algorithm in C and modified it suitably 
based on my suggestions to enumerate the $A$ and $C$ matrices. I would 
like to thank Prof. Bratley for sending me his java implementation of the BM algorithm which showed me 
the merits of efficient programming among many other things. I thank Nicolas Destainville and Naveen Prabhakar for useful comments on an earlier draft of the manuscript.  Last but 
not the least, I would like to thank the High Performance Computing 
Environment at IIT Madras (http://hpce.iitm.ac.in) for providing me 
access to the Leo and Vega superclusters where all the computations were 
carried out.

\appendix

\section{Ferrers Diagrams in a symmetric box}\label{appbox}

Let us consider Ferrers diagrams of i.d. $r$ that fit in a symmetric box 
of size $b$ -- points that lie within the box are such that all their 
coordinates take values in $(0,1,\ldots,b-1)$. Let us call them 
restricted Ferrers Diagrams.  It is easy to see that under the action of 
$S_r$ that permutes the $r$-axes, FD's that fit in a box get mapped to 
FD's that also fit in the same box.  Due to this property, we can 
construct analogs of the various triangles $A/C/D/F$ for restricted FD's 
as well even though the total number of restricted FD's are finite. For 
instance, we have
\begin{equation}
p_d^{\textrm{rest}}(n) =\sum_{r=0}^{d+1} \binom{d+1}{r}\ a^{\textrm{rest}}_{n,r}\ ,
\end{equation}
where $p_d^{\textrm{rest}}(n)$ is the number of FD's with $n$-nodes that 
fit in a symmetric box of size $b^{d+1}$ and 
$A^{\textrm{rest}}=\big(a_{n,r}^{\textrm{rest}}\big)$. Similarly, we can 
define $C^{\textrm{rest}}$. The analog of the $A$-matrix for restricted partitions not necessarily in a symmetric box has appeared in the work of Destainville et. al.\cite{Destainville:1997}.

Let us focus on  partitions that fit into a 
symmetric box of size two and denote the corresponding triangles in obvious notation: 
\begin{center}
$A^{\textrm{box2}}=\big(a^{\textrm{box2}}_{n,r}\big)$, 
$C^{\textrm{box2}}=\big(c^{\textrm{box2}}_{m,x}\big)$, 
$D^{\textrm{box2}}=\big(d^{\textrm{box2}}_{m,x}\big)$ and 
$F^{\textrm{box2}}=\big(f^{\textrm{box2}}_{n,r}\big)$.
\end{center}
 We do not write 
out their relationships as they exactly mirror the corresponding formulae for 
unrestricted partitions.

\begin{mydef}
Let $\mathcal{B}$ be the set of strict Ferrers diagrams that fit in a 
symmetric box of size $2$.
\end{mydef}
We say, somewhat loosely, that a strict skew FD, $\sigma$ of r.d. $x$ is 
in $\mathcal{B}$ if the FD $\mu_x\cup \sigma\in \mathcal{B}$.   The only 
irreducible strict skew Ferrers diagrams at dimensions $2$, $3$ and $4$ in $\mathcal{B}$ are
$\sigma_{2}$, $\sigma_3$ and $\sigma_{4a/b}$ 
respectively. It is also easy to see that $\mathcal{D}\subset \mathcal{B}$.

The matrix $C^{\textrm{box2}}=\big(c^{\textrm{box2}}_{m,x}\big)$ (for 
$m,x\geq0$) has non-zero entries when $m\in[0,2^x-x-1]$ with 
$c^{\textrm{box2}}_{0,0}=1$. Further, 
$c^{\textrm{box2}}_{0,x}=c^{\textrm{box2}}_{x,0}=0$ for $x> 0$ and 
$c^{\textrm{box2}}_{2^x-x-1,x}=1$. Thus, it is a triangle. The maximum 
value of $m$, for fixed $x$, is obtained by considering the FD 
containing all nodes that are in the box. Such an FD has $2^x$ nodes in 
$x$-dimensions and thus the corresponding skew FD has $2^x-x-1$ nodes 
after deleting the nodes that lie in $\mu_x$. Below, we provide the 
first few rows of the matrix, $C^{\textrm{box2}}$, for $m\in [1,6]$.
\begin{equation}
{\scriptsize
C^{\textrm{box2}}=\left(\begin{smallmatrix}
1 \\
0& 0 & 1  \\
 0&0 & 0 & 3 & 3  \\
 0&0 & 0 & 1 & 16 & 30 & 15  \\
 0&0 & 0 & 1 & 15 & 135 & 330 & 315 & 105  \\
 0&0 & 0 & 0 & 18 & 232 & 1581 & 4410 & 5880 & 3780 & 945  \\
 0& 0 & 0 & 0 & 13 & 355 & 4000 & 23709 & 71078 & 116550 & 107100 & 51975 &
   10395  \\
\end{smallmatrix}\right)
}
\end{equation}

Extending the ideas that were used in defining the $D/F$ triangles, we 
look to count \textit{only} those skew FD's that do not contain skew 
FD's in $\mathcal{B}$ are reducible components. Let the matrix 
$\widehat{C}_{m,x}$ denote this reduced $C$-matrix that counts strict 
skew FD's of r.d. $x$ with $m$-nodes. Such FD's necessarily contain at 
least one node of type 2. Then one has the following relation that relates $\widehat{C}$ 
to the $C$:
\begin{equation}
c_{m,x} = \sum_{y=0}^x \sum_{p=0}^{m}\binom{x}{y}\  c^{\textrm{box2}}_{m-p,x-y}\  \widehat{c}_{p,y}\ ,
\end{equation}
with $\widehat{c}_{0,0}\equiv 1$ and $\widehat{c}_{0,x}=0$ for $x>0$. It 
is better to rewrite the above formula as follows:
\begin{equation}\label{chat}
c_{m,x} = \widehat{c}_{m,x} +c^{\textrm{box2}}_{m,x}+\sum_{y=1}^{x-1} \sum_{p=1}^{m-1}\binom{x}{y}\  c^{\textrm{box2}}_{m-p,x-y}\  \widehat{c}_{p,y}\  .
\end{equation}
The first term in the right hand side of the above equation is the 
contribution from skew FD's that do no contain any reducible components 
in $\mathcal{B}$, the second term arise solely from terms that fit into 
a box of size $2$. The last terms runs over terms that contain reducible 
components in $\mathcal{B}$ but do not fit into a box of size $2$. The 
next proposition shows that $\widehat{C}$ is a lower-triangular matrix 
with the $m$-th row containing $m$ non-zero terms.

\begin{prop}\label{chatbound}
$\widehat{c}_{m,x}=0$ when $x>m$ or equivalently when the density $\rho<1$.
\end{prop}
\textbf{Proof:} Since $\mathcal{D}\subset\mathcal{B}$ and all the irreducible strict skew FD's with density less than unity lie in $\mathcal{D}$, the above Proposition follows from Proposition \ref{Fprop}. \mbox{}\hfill  $\Box$

\begin{equation}
\widehat{C}=\left(\begin{smallmatrix}
1 \\
0 &1  \\
 0&1 & 3  \\
0& 1 & 7 & 16  \\
0& 1 & 11 & 57 & 125  \\
0& 1 & 18 & 135 & 602 & 1296 \\
 0&1 & 26 & 293 & 1911 & 7980 & 16807  \\
 0&1 & 38 & 574 & 5242 & 31860 & 127977 & 262144  \\
 0& 1 & 52 & 1089 & 12972 & 106505 & 619872 & 2411416 & 4782969  \\
   \end{smallmatrix} \right)
\end{equation}
We observe that $\widehat{c}_{m,m}= (m+1)^{m-1}$.

We can carry out a similar refinement for strict FD's that contribute to 
the $A$-triangle. One has
\begin{equation}\label{hatF}
\boxed{
a_{n,r} = \widehat{f}_{n,r}+\sum_{s=0}^{r-1} \sum_{p=s+1}^{n-1} \binom{r}s \ a^{\textrm{box2}}_{n-p+1,r-s}\ \widehat{f}_{p,s}
}\ ,
\end{equation} 
with $\widehat{f}_{1,0}\equiv 1$ and $\widehat{f}_{n,0}=0$ for $n>0$. In 
order to interpret the first term, it is better to think of $a_{n,r}$ as 
the number of skew FD's obtained after removing the node at the origin 
of a strict FD. Then, the second term is the contribution from such skew 
FD's that do \textit{not} contain reducible components that fit in a box 
of size two. A second equivalent definition in terms of $m$ is as 
follows:
\begin{equation}
a_{m+r+1,r} = \widehat{f}_{m+r+1,r}+a^{\textrm{box2}}_{m+r+1,r} +\sum_{s=1}^{r-1} \sum_{p=0}^{m} \binom{r}s \ a^{\textrm{box2}}_{m-p+r-s+1,r-s}\ \widehat{f}_{p+s+1,s}\ . 
\end{equation}

It is easy to see there is a bijective map that relates skew FD's that 
contribute to $\widehat{c}_{m,x}$ and those that contribute to 
$\widehat{a}_{m+x+1,x}$. The bijection follows by observing that if 
$\sigma$ is a strict skew FD with $m$ nodes and r.d. $x$, there is a 
unique FD (with i.d. and r.d. equal to $x$) obtained by adding the nodes in $\mu_x$.  Thus,
$$
\widehat{f}_{m+x+1,x}=\widehat{c}_{m,x}\ .
$$

It is easy to see using Proposition \ref{chatbound} that for $\widehat{f}_{n,r}=0$ when $r> n/2$. We 
define the matrix $\widehat{F}=\big(\widehat{f}_{n,r})$ for 
$n=1,2,\ldots$ and $r=0,1,2,\ldots$. Further, we observe that 
$\widehat{f}_{2x+1,x}=c_{x,x}=(x+1)^{x-1}$. Below we reproduce the first 
eleven rows of the $\widehat{F}$-matrix. We reproduce the $F$-matrix alongside for comparison. The first instance where they differ is when $n=8$ and $r=3$ -- this is precisely where the node $(1,1,1)^T$ that is not in $\mathcal{D}$ but present in $\mathcal{B}$ appears. As we go to higher values of $n$,   an entry in $\widehat{F}$ will be generically smaller than the corresponding entry in $F$. Further, entries in both matrices will agree when the density is in the range $[1,4/3)$ -- this is because a node present in $\mathcal{B}$ but not in $\mathcal{D}$ first appears at density $4/3$.
\begin{equation}
\widehat{F} =\left( \begin{smallmatrix}
1\\
 0  \\
0& 1  \\
0& 1  \\
0& 1 & 3  \\
0& 1 & 7  \\
0& 1 & 11 & 16  \\
0& 1 & 18 & 57\\
0& 1 & 26 & 135 & 125  \\
0& 1 & 38 & 293 & 602 \\
0& 1 & 52 & 574 & 1911 & 1296  \\
 \end{smallmatrix}\right)\quad,\quad
F=\left(\begin{smallmatrix}
1\\
0  \\
 0&1  \\
0& 1  \\
 0&1 & 3  \\
 0&1 & 7  \\
 0&1 & 11 & 16  \\
 0&1 & 18 & 58  \\
 0& 1 & 26 & 135 & 125  \\
 0& 1 & 38 & 293 & 618 \\
 0& 1 & 52 & 574 & 1927 & 1296 \\
   \end{smallmatrix}\right)
\ .
 \end{equation}
The second row has only vanishing entries. That is because the only 
strict FD with two nodes fits in a box of size two. So the first 
non-vanishing contribution appears at $n=3,r=1$ if we ignore the 
$n=1,r=0$ term that is more less part of the definition.
 
We can now revisit the problem of enumerating partitions of $n$ in any 
dimension. We see that we need to enumerate the first $n$ rows of the 
$\widehat{F}$ matrix and $A^{\textrm{box2}}$ in order to obtain row $n$ 
of the $A$-matrix. However, from Eq. \eqref{hatF} we see that it is 
sufficient to determine only the first $[n/2]$ elements of row $n$ as 
that completely determines row $n$ of $\widehat{F}$. However, this reduction
is accompanied by the need to evaluate $A^{\textrm{box2}}$  which is yet another computation.
Hence, we preferred to work with the $F$-matrix. However, one should be open to 
using the $\widehat{F}$-matrix if one has an algorithm to directly compute it. Then, the additional effort to compute $A^{\textrm{box2}}$ might be worth it.

\subsection{The box transform}

Define the following generating function for the $A$-matrix
\begin{equation}\label{boxtransform}
A(q,t)= \sum_{m=0}^\infty \sum_{r=0}^\infty a_{m+r+1,r}\ \frac{q^m t^r}{r!}\ ,
\end{equation}
along with similar definitions for $A^{\textrm{box2}}(q,t)$ and 
$\widehat{F}(q,t)$. Then, Eq. \eqref{hatF} implies that the generating 
functions have a simple relation. One has
\begin{equation}
A(q,t) = A^{\textrm{box2}}(q,t) \times \widehat{F}(q,t)\ .
\end{equation} 
It is due to this property that we refer to Eq. \eqref{boxtransform} as 
the \textit{box transform}. Similarly, one defines
\begin{equation}
C(q,t)= \sum_{m=0}^\infty \sum_{r=0}^\infty c_{m,r}\  \frac{q^m t^r}{r!}\ ,
\end{equation}
along with similar definitions for $C^{\textrm{box2}}(q,t)$ and 
$\widehat{C}(q,t)$. Again, one has
\begin{equation}
C(q,t) = C^{\textrm{box2}}(q,t) \times \widehat{C}(q,t)\ .
\end{equation} 

There is an obvious extension to our considerations by replacing the 
symmetric box of size two by one of size $b$. Again, relations of the 
kind that we considered between FD's that fit in the box and those that 
don't appear. For instance, one has
\begin{equation}
A(q,t) = A^{\textrm{box}b}(q,t) \times \widehat{F}(q,t)\ ,
\end{equation} 
where $\widehat{A}(q,t)$ is the generating function of FD's that don't 
fit into a box of size $b$ and do not have reducible parts that fit into 
the box.

\bibliography{partitions}

\clearpage
\begin{sidewaystable}
\footnotesize
\begin{tabular}{c|rrrrrrrrrr} \hline
$n$ & $d=1$ &  $d=2$ &  $d=3$ &  $d=4$ &  $d=5$ &  $d=6$ &  $d=7$ &  $d=8$ &  $d=9$ &  $d=10$ \\ \cline{2-11}
1&  1 & 1 & 1 & 1 & 1 & 1 & 1 & 1 & 1 & 1 \\
 2 & 2 & 3 & 4 & 5 & 6 & 7 & 8 & 9 & 10 & 11 \\
3&  3 & 6 & 10 & 15 & 21 & 28 & 36 & 45 & 55 & 66 \\
4&  5 & 13 & 26 & 45 & 71 & 105 & 148 & 201 & 265 & 341 \\
5& 7 & 24 & 59 & 120 & 216 & 357 & 554 & 819 & 1165 & 1606 \\
 6&11 & 48 & 140 & 326 & 657 & 1197 & 2024 & 3231 & 4927 & 7238 \\
 7&15 & 86 & 307 & 835 & 1907 & 3857 & 7134 & 12321 & 20155 & 31548 \\
8&  22 & 160 & 684 & 2145 & 5507 & 12300 & 24796 & 46209 & 80920 & 134728 \\
 9& 30 & 282 & 1464 & 5345 & 15522 & 38430 & 84625 & 170370 & 319555 & 565983 \\
 10 & 42 & 500 & 3122 & 13220 & 43352 & 118874 & 285784 & 621316 & 1247780 & 2350183 \\
11&  56 & 859 & 6500 & 32068 & 119140 & 362670 & 953430 & 2240838 & 4821737 & 9661465 \\
12 & 77 & 1479 & 13426 & 76965 & 323946 & 1095430 & 3151332 & 8011584 & 18478640 & 39401792 \\
 13 & 101 & 2485 & 27248 & 181975 & 869476 & 3271751 & 10314257 & 28395213 & 70261505 & 159527302 \\
14 & 135 & 4167 & 54804 & 425490 & 2308071 & 9673993 & 33457972 & 99845553 & 265266530 & 641733862 \\
 15 &176 & 6879 & 108802 & 982615 & 6056581 & 28310881 & 107557792 & 348333411 & 994606250 & 2565774277 \\
 16 &231 & 11297 & 214071 & 2245444 & 15724170 & 82033609 & 342732670 & 1205925033 & 3704360354 & 10198601886 \\
17 & 297 & 18334 & 416849 & 5077090 & 40393693 & 235359901 & 1082509680 & 4142850423 & 13705110470 & 40305279454 \\
 18 & 385 & 29601 & 805124 & 11371250 & 102736274 & 668779076 & 3389190112 & 14122999548 & 50367905030 & 158376907546 \\
 19 &490 & 47330 & 1541637 & 25235790 & 258790004 & 1882412994 & 10518508294 & 47772540002 & 183864216415 & 618742851276 \\
 20 &627 & 75278 & 2930329 & 55536870 & 645968054 & 5249817573 & 32361863632 & 160336300356 & 666612686420 & 2403142436321 \\
21 & 792 & 118794 & 5528733 & 121250185 & 1598460229 & 14510628853 & 98711666690 & 533909133114 & 2400146830007 \\
 22 & 1002 &  186475 & 10362312 & 262769080 & 3923114261& 39762851345 &  298546248070 &
  1763901729589 & 8581152930795 \\
23 & 1255 & 290783 & 19295226 & 565502405 & 9554122089 & 108058883583 &  895425789360\\ 
\end{tabular}
\caption{$d\leq 10$-dimensional partitions of $n\leq 23$ as determined by direct enumeration using Knuth's algorithm. This provides an independent cross-check of the entries in the first $11$ columns of the $A$-matrix.}\label{Knuthresults}
\end{sidewaystable}
\newpage

\begin{sidewaystable}
%\scriptsize
\begin{tabular}{r|r@{~}r@{~}r@{~}r@{~}r@{~}r@{~}r@{~}r@{~}r@{~}r@{~}r@{~}r@{~}r@{~}r@{~}r@{~}r@{~}r@{~}r@{~}r@{~}r@{~}r@{~}}
%\begin{array}{cccccccccccccccccccccc}
%$n$ & $0$ &$1$ &$2$ &$3$ &$4$ &$5$ &$r=6$ &$r=7$ &$r=8$ &$r=9$ &$r=10$ &$r=11$ &$r=12$ &$r=13$ &$r=14$ &$r=15$ &$r=16$ &$r=17$ &$r=18$ &$r=19$ \\
$n\backslash r$\!\!   & $0$ & $1$ &$2$ &$3$ &$4$ &$5$ &$6$ &$7$ &$8$ &$9$ &$10$ &$11$ 
%&$12$ &$13$ &$14$ &$15$ &$16$ &$17$ &$18$ 
\\ \hline
1 & 1 \\
2& 0&  1 \\
 %& 0 & 0 & 0 & 0 & 0 & 0 & 0 & 0 & 0 & 0 & 0 & 0 & 0 & 0 & 0 & 0 & 0 & 0 & 0 & 0  \\
3& 0&  1 & 1 \\ 
 %& 0 & 0 & 0 & 0 & 0 & 0 & 0 & 0 & 0 & 0 & 0 & 0 & 0 & 0 & 0 & 0 & 0 & 0 & 0  \\
4& 0&   1 & 3 & 1 \\
 %& 0 & 0 & 0 & 0 & 0 & 0 & 0 & 0 & 0 & 0 & 0 & 0 & 0 & 0 & 0 & 0 & 0 & 0  \\
 5&0&   1 & 5 & 6 & 1 \\
 %& 0 & 0 & 0 & 0 & 0 & 0 & 0 & 0 & 0 & 0 & 0 & 0 & 0 & 0 & 0 & 0 & 0  \\
6& 0&  1 & 9 & 18 & 10 & 1 \\
 %& 0 & 0 & 0 & 0 & 0 & 0 & 0 & 0 & 0 & 0 & 0 & 0 & 0 & 0 & 0 & 0 \\
7&0&   1 & 13 & 44 & 49 & 15 & 1 \\
 %& 0 & 0 & 0 & 0 & 0 & 0 & 0 & 0 & 0 & 0 & 0 & 0 & 0 & 0 & 0  \\
8&0&   1 & 20 & 97 & 172 & 110 & 21 & 1 \\
 %& 0 & 0 & 0 & 0 & 0 & 0 & 0 & 0 & 0 & 0 & 0 & 0 & 0 & 0\\
9& 0&   1 & 28 & 195 & 512 & 550 & 216 & 28 & 1 \\
 %& 0 & 0 & 0 & 0 & 0 & 0 & 0 & 0 & 0 & 0 & 0 & 0 & 0  \\
10&0&   1 & 40 & 377 & 1370 & 2195 & 1486 & 385 & 36 & 1 \\
 %& 0 & 0 & 0 & 0 & 0 & 0 & 0 & 0 & 0 & 0 & 0 & 0  \\
11& 0&  1 & 54 & 694 & 3396 & 7603 & 7886 & 3514 & 638 & 45 & 1 \\
 %& 0 & 0 & 0 & 0 & 0 & 0 & 0 & 0 & 0 & 0 & 0  \\
12& 0&  1 & 75 & 1251 & 7968 & 23860 & 35115 & 24318 & 7484 & 999 & 55 & 1 \\
 %& 0 & 0 & 0 & 0 & 0 & 0 & 0 & 0 & 0 & 0  \\
13&0&   1 & 99 & 2185 & 17910 & 69580 & 138155 & 138075 & 65997 & 14667 & 1495 & 66 \\
 %& 1 & 0 & 0 & 0 & 0 & 0 & 0 & 0 & 0 & 0  \\
 14 &0&  1 & 133 & 3765 & 38942 & 191795 & 495870 & 677663 & 471276 & 161202 & 26875 & 2156  \\
 %& 78 & 1& 0 & 0 & 0 & 0 & 0 & 0 & 0 & 0  \\
 15&0&   1 & 174 & 6354 & 82338 & 505640 & 1657975 & 2978735 & 2864408 & 1424142 & 360940 & 46596  \\
 %& 3015 & 91 & 1& 0 & 0 & 0 & 0 & 0 & 0 & 0  \\
16&0&    1 & 229 & 10607 & 170265 & 1285754 & 5240090 & 12016809 & 15354492 & 10604286 & 3880561 & 751696  \\
 %& 77144 & 4108 & 105 & 1& 0 & 0 & 0 & 0 & 0 & 0  \\
 17 &0&  1 & 295 & 17446 & 345291 & 3173220 & 15821657 & 45268685 & 74497870 & 68869266 & 34954135 & 9685709 \\
 %& 1472923 & 122824 & 5474 & 120 & 1 & 0 & 0 & 0 & 0 & 0 \\
18&0&   1 & 383 & 28449 & 689026 & 7637795 & 45999383 & 161270025 & 333494972 & 400292769 & 272579245 & 104184949  \\ 
% & 22435260 & 2740322 & 189112 & 7155 & 136 & 1& 0 & 0 & 0 & 0 \\
19&0&   1 & 488 & 45863 & 1355253 & 17996010 & 129560563 & 548523528 & 1397398036 & 2123894171 & 1886698315 & 965585764 \\
   % & 284668913 & 48745112 & 4876067 & 282850 & 9196 & 153 & 1& 0 & 0 & 0  \\
20&0&    1 & 625 & 73400 & 2632975 & 41631740 & 355205608 & 1794375520 & 5541288850 & 10446368715 & 11819801575 & 7897875909  \\
   %& 3102741280 & 721400004 & 100211944 &8347235 & 412456 & 11645 & 171&1 & 0 & 0 \\
21 &0&   1 & 790 & 116421 & 5058305 & 94786545 & 951526108 & 5678296645 & 20973892932 & 48206965521 & 68073453307 & 58101011914  \\
%& 29731926984 & 9153619020 & 1712287591 &196337499 & 13813806 & 588149 & 14553& 1    & 0 
22 &0&   1& 1000  &183472 & 9622420& 212812255 &2498219985 &17463026868 & 76290515426 &210725428060&
  364964576905& 390349624764 \\
%& 255102172260 & 101855999739 & 25054962934 & 3838014650 & 368763408 & 22187703 & 822189 & 17974 & 210 \\
  23 &0&   1 &1253 &287021 &18139620 &471921560 &6444739208 &52390397612
   &268136421612 &879260678868 &1840128105650 &2425318710876 \\
%   &1989135984936 &1011712387244 &320976466643 &64203399495 &8179916108 &667249575 &34704447 \\
   %&1129132 &21965 \\
\end{tabular}
\caption{The first $12$ columns and $23$ rows of the triangle $A$.  The other $11$ columns can be obtained using the ten rows of the $C$-matrix given below. Thus, one can determine partitions of positive integers $\leq 23$ from it.}\label{Amatrix}
\bigskip
%\end{sidewaystable}
%\begin{sidewaystable}
\tiny
\begin{tabular}{r|r@{~}r@{~}r@{~}r@{~}r@{~}r@{~}r@{~}r@{~}r@{~}r@{~}r@{~}r@{~}r@{~}r@{~}r@{~}r@{~}r@{~}r@{~}r@{~}r@{~}r@{~}}
$m\backslash x$\!\! & $0$ & $1$ &$2$ &$3$ &$4$ &$5$ &$6$ &$7$ &$8$ &$9$ &$10$ &$11$ &$12$ &$13$ &$14$ &$15$ &$16$ &$17$ &$18$ &$19$  \\ \hline
0 & 1 \\
 1&0 & 1 & 1  \\
 2&0 & 1 & 3 & 6 & 3  \\
 3&0 &1 & 7 & 20 & 46 & 45 & 15  \\
 4&0 &1 & 11 & 61 & 198 & 480 & 645 & 420 & 105  \\
 5&0 &1 & 18 & 138 & 706 & 2508 & 6441 & 10395 & 9660 & 4725 & 945  \\
 6&0 &1 & 26 & 296 & 2052 & 10375 & 38809 & 105392 & 192668 & 224595 & 159075
   & 62370 & 10395 \\
 7&0 &1 & 38 & 577 & 5428 & 36285 & 184624 & 713402 & 2032500 & 4080195 &
   5580855 & 5051970 & 2889810 & 945945 & 135135  \\
 8&0 &1 & 52 & 1092 & 13226 & 114220 & 751639 & 3854487 & 15231326 & 45159822
   & 97613505 & 150613155 & 162889650 & 120270150 & 57702645 & 16216200 &
   2027025  \\
 9&0 &1 & 73 & 1963 & 30648 & 332035 & 2747799 & 17918432 & 92357844 &
   370929320 & 1136808010 & 2609559315 & 4427605050 & 5488733250 &
   4892112225 & 3047969925 & 1259458200 & 310134825 & 34459425  \\
 10&0 & 1 & 97 & 3471 & 67868 & 910729 & 9268382 & 74767133 & 483797592 &
   2498431224 & 10155656364 & 31998207087 & 77214286182 & 141528086700 &
   195617897475 & 201837365730 & 152796603960 & 82323566325&
   29876321475 & 6547290750  \\
\end{tabular}
\caption{The second triangle -- the first ten rows  and nineteen columns of the $C$-matrix. We have only shown non-zero entries. 
}\label{Cmatrix}
\end{sidewaystable}
\begin{sidewaystable}
\footnotesize
\begin{tabular}{r|r@{~}r@{~}r@{~}r@{~}r@{~}r@{~}r@{~}r@{~}r@{~}r@{~}r@{~}r@{~}r@{~}r@{~}r@{~}r@{~}r@{~}r@{~}r@{~}r@{~}r@{~}}
$n\backslash x$\!\! & $0$ & $1$ &$2$ &$3$ &$4$ &$5$ &$6$ &$7$ &$8$ &$9$ &$10$ &$11$ &$12$   \\ \hline
1 & 1 \\
2 & 0 \\
3  & 0& 1  \\
 4 & 0 &1  \\
 5  & 0&1 & 3  \\
 6 & 0&1 & 7  \\
 7 & 0& 1 & 11 & 16  \\
 8 & 0& 1 & 18 & 58  \\
 9 & 0&1 & 26 & 135 & 125  \\
 10  & 0&1 & 38 & 293 & 618  \\
11 & 0 & 1 & 52 & 574 & 1927 & 1296 \\
 12  & 0&1 & 73 & 1089 & 5256 & 8220  \\
 13  & 0&1 & 97 & 1960 & 12982 & 32380 & 16807  \\
 14  & 0&1 & 131 & 3468 & 30320 & 107270 & 131897 \\
 15  & 0&1 & 172 & 5955 & 67414 & 319530 & 633442 & 262144  \\
 16 & 0 &1 & 227 & 10085 & 145045 & 888983 & 2490187 & 2483096  \\
 17 & 0 &1 & 293 & 16759 & 303101 & 2346515 & 8710068 & 14200018 & 4782969  \\
 18  & 0&1 & 381 & 27564 & 619564 & 5952280 & 28205459 & 65151254 & 53672292  \\
 19 & 0 &1 & 486 & 44714 & 1241845 & 14617100 & 86238209 & 263040064 & 359302890
   & 100000000  \\
 20 & 0 &1 & 623 & 71936 & 2450043 & 34962755 & 252190709 & 975528302 &
   1899997612 & 1309707840  \\
 21 & 0 &1 & 788 & 114546 & 4765327 & 81792100 & 711409264 & 3398678150 &
   8749699709 & 10128660960 & 2357947691  \\
 22  & 0&1 & 998 & 181102 & 9157550 & 187791450 & 1948153500 & 11278286646 &
   36739765288 & 61114773760 & 35600917115  \\
 23 & 0 &1 & 1251 & 284021 & 17406714 & 424233500 & 5203415684 & 35979941641 &
   144179174632 & 318163092360 & 314636749085 & 61917364224  \\
 24 & 0 &1 & 1571 & 442713 & 32771292 & 944990470 & 13605818265 & 111092074842 &
   536798419714 & 1499829016296 & 2148096711540 & 1066426694784 \\
 25  & 0&1 & 1954 & 685443 & 61158328 & 2079070155 & 34930133300 & 333670251012 &
   1915118952548 & 6574308285588 & 12551603978445 & 10672681371264 &
   1792160394037 
\end{tabular}
\caption{The $F$-matrix as determined using data up to $a_{25,12}$. This determines partitions of all integers $\leq 25$. We have only shown non-zero entries.}\label{Ftable}
\end{sidewaystable}

\begin{sidewaystable}
\footnotesize
\begin{tabular}{r|r@{~}r@{~}r@{~}r@{~}r@{~}r@{~}r@{~}r@{~}r@{~}r@{~}r@{~}r@{~}r@{~}r@{~}r@{~}r@{~}r@{~}r@{~}r@{~}r@{~}r@{~}}
$m\backslash x$\!\!  & $0$ & $1$ &$2$ &$3$ &$4$ &$5$ &$6$ &$7$ &$8$ &$9$ &$10$ &$11$ &$12$ &$13$ &$14$ &$15$  \\ \hline
0 & 1 \\
1&0 &1  \\
2&0 & 1 & 3 & 3  \\
3&0 & 1 & 7 & 17 & 28  \\
4&0 & 1 & 11 & 58 & 156 & 295 & 90  \\
5&0 & 1 & 18 & 135 & 640 & 1913 & 3786 & 2310  \\
6&0 & 1 & 26 & 293 & 1944 & 9010 & 28714 & 59024 & 50960 & 7560 \\
7&0 & 1 & 38 & 574 & 5272 & 33340 & 154654 & 509912 & 1089488 & 1158192 &
   378000  \\
8&0 & 1 & 52 & 1089 & 12998 & 108465 & 671389 & 3123477 & 10485214 & 23226165
   & 28428750 & 14206500 & 1247400 \\
9&0 & 1 & 73 & 1960 & 30336 & 321130 & 2551119 & 15580292 & 72440912 &
   245511503 & 561332710 & 762518790 & 501491760 & 102702600  \\
10&0 & 1 & 97 & 3468 & 67430 & 891114 & 8811002 & 67908953 & 409620720 &
   1895816757 & 6456110604 & 15166699372 & 22350118032 & 17852174340 &
   5864859000 & 340540200  \\
%11 & 1 & 131 & 5955 & 145061 & 2349415 & 28374212 & 268506497 & 2021849016 &
%   12029683461 & 55244054100 & 188494710806 & \text{DD}(11,12) &
%   \text{DD}(11,13) & \text{DD}(11,14) & \text{DD}(11,15) &
%   \text{DD}(11,16)  \\
% 12 &1 & 172 & 10085 & 303117 & 5956110 & 86523342 & 986425727 & 9033780362 &
%   66589017808 & 390098243490 & \text{DD}(12,11) & \text{DD}(12,12) &
%   \text{DD}(12,13) & \text{DD}(12,14) & \text{DD}(12,15) &
%   \text{DD}(12,16) & \text{DD}(12,17) & 138940401600 \\
   \end{tabular}
\caption{The first eleven rows of the $D$-matrix. We have only shown non-zero entries.
}\label{Dmatrix}
\end{sidewaystable}
\begin{sidewaystable}
\scriptsize
\begin{tabular}{r|r@{~}r@{~}r@{~}r@{~}r@{~}r@{~}r@{~}r@{~}r@{~}r@{~}r@{~}r@{~}r@{~}r@{~}r@{~}r@{~}r@{~}r@{~}r@{~}r@{~}r@{~}}
$z\backslash y$\!\! & $0$ & $1$ &$2$ &$3$ &$4$ &$5$ &$6$ &$7$ &$8$ &$9$ &$10$   \\ \hline
0&1 & 0 & 0 & 0 & 0 & 0 & 0 & 0 & 0 & 0 & 0 \\
 1 &1 & 0 & 0 & 0 & 0 & 0 & 0 & 0 & 0 & 0 & 0 \\
 2 &$-1$ & 0 & 0 & 0 & 0 & 0 & 0 & 0 & 0 & 0 & 0 \\
 3 &$-5$ & 6 & 0 & 0 & 0 & 0 & 0 & 0 & 0 & 0 & 0 \\
4& 1 & 74 & 0 & 0 & 0 & 0 & 0 & 0 & 0 & 0 & 0 \\
5& 41 & 252 & $-40$ & 0 & 0 & 0 & 0 & 0 & 0 & 0 & 0 \\
 6&31 & $-540$ & $-676$ & 0 & 0 & 0 & 0 & 0 & 0 & 0 & 0 \\
 7&$-461$ & $-6470$ & 1180 & 656 & 0 & 0 & 0 & 0 & 0 & 0 & 0 \\
 8& $-895$ & $-4074$ & 69020 & 10864 & 0 & 0 & 0 & 0 & 0 & 0 & 0 \\
 9&6481 & 138264 & 403620 & $-39016$ & $-9216$ & 0 & 0 & 0 & 0 & 0 & 0 \\
 10 &22591 & 376872 & $-961240$ & $-1628984$ & $-191456$ & 0 & 0 & 0 & 0 & 0 & 0 \\
 11 &$-107029$ &$ -2922930$ & $-21162456$ & $-3687040$ & 463680 & 195840 & 0 & 0 & 0 & 0 & 0 \\
 12 & $\beta_{12,0}$ &$ -15867390$ & $-40350840$ & 168546560 & 40336016 & 7455104 & 0 & 0 & 0 & 0 & 0 \\
13 & $\beta_{13,0}$ & $\beta_{13,1}$ & 758778240 & 1656046448 & 110435472 & 73922176 & $-6297600$ & 0 & 0 & 0 & 0 \\
 14 & $\beta_{14,0}$ & $\beta_{14,1}$ & $\beta_{14,2}$ & $-1927766192$ & $-5730022032$ & $-552798336$ & $-382393600$ & 0 & 0 & 0 & 0 \\
 15 & $\beta_{15,0}$ & $\beta_{15,1}$ &  $\beta_{15,2}$ &  $\beta_{15,3}$ & $-44646818832$ & $-10585577760$ & $-7549384960$ & 278906880 & 0 & 0 & 0 \\
16 & $\beta_{16,0}$& $\beta_{16,1}$ & $\beta_{16,2}$ & $\beta_{16,3}$ & $\beta_{16,4}$ & 75450085920 & $-14753227264$ & 22686050304 & 0 & 0 & 0 \\
17 & $\beta_{17,0}$&  $\beta_{17,1}$ &  $\beta_{17,2}$&  $\beta_{17,3}$ &  $\beta_{17,4}$ & $\beta_{17,5}$ & 1603141023616 & 607200778752 &$-14729379840$ & 0 & 0 \\
18 & $\beta_{18,0}$ & $\beta_{18,1}$ & $\beta_{18,2}$ & $\beta_{18,3}$& $\beta_{18,4}$ & $\beta_{18,5}$ & $\beta_{18,6}$ & 2727351931392 & $-1449282760704$ & 0 & 0 \\
19 &  $\beta_{19,0}$ &$\beta_{19,1}$& $\beta_{19,2}$ & $\beta_{19,3}$ &$\beta_{19,4}$ & $\beta_{19,5}$ & $\beta_{19,6}$ & $\beta_{19,7}$ & $-47662776674304$ & 873791815680 & 0 \\
20 & $\beta_{20,0}$ & $\beta_{20,1}$ & $\beta_{20,2}$ & $\beta_{20,3}$ & $\beta_{20,4}$ & $\beta_{20,5}$& $\beta_{20,6}$ &$\beta_{20,7}$ & $\beta_{20,8}$ & 101710939668480 & 0 \\
21 & $\beta_{21,0}$ & $\beta_{21,1}$ & $\beta_{21,2}$& $\beta_{21,3}$ & $\beta_{21,4}$ & $\beta_{21,5}$& $\beta_{21,6}$ & $\beta_{21,7}$ & $\beta_{21,8}$ & $\beta_{21,9}$ & $-58358690611200$
  \end{tabular}
\caption{The $\beta$-triangle to the extent that we have determined it. The first column is consistent with the Meeussen conjecture. The $\beta_{m+x,x}$ for $x\in[0,m-1]$ completely determine the degree $2m$ polynomial $g_{m}(r)$. Thus, we have determined all polynomials for $m\leq 11$ albeit assuming the existence of the $\beta$-matrix which is conjectural. The polynomials obtained this way agrees with the ones determined by the $C/D/F$ matrices.
}\label{Betamatrix}
\end{sidewaystable}
\begin{sidewaystable}
\tiny
\begin{tabular}{r|r@{~}r@{~}r@{~}r@{~}r@{~}r@{~}r@{~}r@{~}r@{~}r@{~}r@{~}r@{~}r@{~}r@{~}r@{~}r@{~}r@{~}r@{~}r@{~}r@{~}r@{~}}
$m\backslash x$\!\! & $0$&  $1$ &$2$ &$3$ &$4$ &$5$ &$6$ &$7$ &$8$ &$9$ &$10$ &$11$ &$12$ &$13$ &$14$ &$15$ &$16$ & $17$ & $18$ & $19$ & $20$ \\ \hline
0 & 1 \\
1&0 &   0 & 1  \\
 2 &0 & 0 & 0 & 3 & 3  \\
 3 &0 & 0 & 0 & 1 & 16 & 30 & 15  \\
 4 &0 & 0 & 0 & 1 & 15 & 135 & 330 & 315 & 105 \\
 5 &0 & 0 & 0 & 0 & 18 & 232 & 1581 & 4410 & 5880 & 3780 & 945  \\
 6&0 &  0 & 0 & 0 & 13 & 355 & 4000 & 23709 & 71078 & 116550 & 107100 & 51975 &
   10395  \\
 7&0 &  0 & 0 & 0 & 10 & 450 & 8075 & 78725 & 431460 & 1353240 & 2552130 &
   2962575 & 2079000 & 810810 & 135135  \\
 8&0 &  0 & 0 & 0 & 6 & 530 & 14065 & 204540 & 1767045 & 9207945 & 29811330 &
   62179425 & 85270185 & 76621545 & 43513470 & 14189175 & 2027025  \\
 9& 0 & 0 & 0 & 0 & 4 & 580 & 22315 & 456400 & 5704580 & 44793784 & 225211165 &
   746795775 & 1680747090 & 2612970360 & 2812925115 & 2062160100 &
   983782800 & 275675400 & 34459425  \\
 10 &0 & 0 & 0 & 0 & 1 & 611 & 33177 & 918981 & 15738310 & 174240318 & 1268511894
   & 6207749790 & 20975922462 & 50107517460 & 85928953110 & 106306245045
   & 94166932860 & 58305347100 & 23983759800 & 5892561675 & 654729075
    \end{tabular}
\caption{The first eleven rows of the $C^{\textrm{box2}}$-triangle. We have only shown non-zero entries.
}\label{Bhatmatrix}
\bigskip
\centering
\scriptsize
\begin{tabular}{r|r@{~}r@{~}r@{~}r@{~}r@{~}r@{~}r@{~}r@{~}r@{~}r@{~}r@{~}r@{~}r@{~}r@{~}r@{~}r@{~}r@{~}r@{~}r@{~}r@{~}r@{~}}
$n\backslash r$\!\! &  $1$ &$2$ &$3$ &$4$ &$5$ &$6$ &$7$ &$8$ &$9$ &$10$ &$11$   \\ \hline
1 & 1\\
2 &0  \\
3 &0  \\
4 & 0  \\
5& 0  \\
6&  0  \\
 7 &0  \\
8 & 0 & 0 & 0 & 1  \\
 9 &0 \\
10& 0 & 0 & 0 & 0 & 12 \\
11& 0 & 0 & 0 & 0 & 12  \\
 12 &0&0 & 0 & 0 & 10 & 150 \\
 13& 0 &0 & 0 & 0 & 6 & 330  \\
 14& 0 &0 & 0 & 0 & 4 & 485 & 2160  \\
 15& 0 &0 & 0 & 0 & 1 & 570 & 7750  \\
 16& 0 &0 & 0 & 0 & 1 & 610 & 17280 & 36015  \\
 17& 0 &0 & 0 & 0 & 0 & 600 & 30120 & 185430 \\
 18& 0 &0 & 0 & 0 & 0 & 580 & 45720 & 574280 & 688128  \\
 19& 0 &0 & 0 & 0 & 0 & 530 & 63870 & 1364195 & 4727520  \\
 20 & 0&0 & 0 & 0 & 0 & 470 & 85325 & 2751875 & 19192880 & 14880348  \\
 21& 0 &0 & 0 & 0 & 0 & 387 & 110625 & 4994640 & 59080000 & 130094748  \\
 22& 0 &0 & 0 & 0 & 0 & 310 & 140322 & 8480885 & 152220320 & 664737850 &
   360000000 
   \\
23& 0 & 0 & 0 & 0 & 0 & 215 & 174380 & 13808620 & 346973284 & 2557358244 &
   3873139200 \\
24 & 0&  0 & 0 & 0 & 0 & 155 & 212815 & 21879725 & 726316080 & 8167776498 &
   24169328400 & 9646149645 
\end{tabular}
\caption{The first $24$ rows of the $F^{\textrm{box2}}$-triangle. We have only shown non-zero entries except for rows which have all zeros where we shown a zero in the first column.
}\label{Fboxtable}
\end{sidewaystable}
\end{document}